\documentclass[11pt]{amsart}
\usepackage[margin=1in]{geometry}
\usepackage{amsmath,amssymb,mathtools}
\usepackage[hidelinks]{hyperref}
\hypersetup{pdftitle={An Obukhov-Corrsin-Yaglom law for a kinematic model of turbulence},pdfauthor={Camillo De Lellis}}
\newtheorem{theorem}{Theorem}[section]
\newtheorem{proposition}[theorem]{Proposition}
\newtheorem{lemma}[theorem]{Lemma}
\newtheorem{corollary}[theorem]{Corollary}
\theoremstyle{definition}
\newtheorem{definition}[theorem]{Definition}
\theoremstyle{remark}
\newtheorem{remark}[theorem]{Remark}
\newcommand{\T}{\mathbb T}
\newcommand{\E}{\mathbb E}
\newcommand{\dd}{\,\mathrm d}
\newcommand{\Id}{\mathrm{Id}}
\newcommand{\norm}[1]{\lVert #1\rVert}
\newcommand{\ind}{\mathbf 1}
\DeclareMathOperator{\diver}{div}
\title[Obukhov-Corrsin-Yaglom Law]{An Obukhov-Corrsin-Yaglom law for a kinematic model of turbulence}
\author{Camillo De Lellis}

\begin{document}
\begin{abstract}
We prove an Obukhov-Corrsin-Yaglom-type law for passive scalars
advected by a family of random $C^\alpha$ velocity field proposed in \cite{JDEMS} as a kinematic model of turbulent fluids. Up to a technical assumption on one parameter, the expectation of the squared $L^2_t H^\gamma_x$ norm of solutions to the passive scalar equation is bounded by the squared $L^2$ norm of the deterministic initial datum, uniformly in the viscosity, whenever $\alpha + 2\gamma < 1$. The proof is based on a new insight in the theory developed by \cite{FM} to prove Richardson's dispersion and anomalous dissipation for the same model. 
\end{abstract}
\maketitle

\section{Introduction}

The Obukhov--Corrsin theory and Yaglom's law
\cite{Obukhov,Corrsin,Yaglom} concern the statistics of passive scalar
turbulence. At the level of regularity, they suggest that a random turbulent
velocity field of a given spatial regularity induces a corresponding
regularity of the scalars it transports, uniformly as diffusion
tends to zero. We refer to
\cite{Warhaft,FGV,ShraimanSiggia,Sreenivasan,AVSurvey} for surveys and further
references. For mathematical developments on mixing and deterministic
anomalous dissipation see
\cite{CKRZ,ACMMixing,ACMLoss,DEIJ,CCS,AV,BSWEuler,BSWRegularity,EL,DE,JSAutonomous,HCRUniversal,HCRIntermittent},
and \cite{BDL,BCCDL,JSMeasure} for related results for fluid equations.
For stochastic models we refer to
\cite{Kraichnan,GK,BGK,LJR,BBPSLagrangian,BBPSBatchelor,HPZZ,GGM,DGP,RowanDiffusion,RowanSpectrum}.

A comparable regularity theory for rough random velocity fields remains much less developed, one notable exception is the recent groundbreaking paper \cite{CHR}. This note contributes
to the latter topic for the kinematic model introduced in the author's joint
work \cite{JDEMS}: building on the recent breakthrough of Fefferman and
Mayboroda \cite{FM}, we prove a bound in expectation for the squared
$L^2_tH^\gamma_x$ norm of passively advected scalars, uniformly in the
diffusivity, whenever the velocity is $C^\alpha$ and
$\alpha+2\gamma<1$ (our statement requires, however, some suitable restrictions on the parameters of the model). The bound
holds for every deterministic $L^2$ initial datum. 

\subsection{The kinematic model} We first detail the kinematic model of \cite{JDEMS}, which constructs certain probability measures on the space of time-dependent divergence-free velocity fields. In fact the model has a lot of flexibility, but for the purpose of this paper we will make some precise choices, following \cite[Section 1.2]{FM}.
 
The velocity is a sum of spatial wavelets whose centers are transported
by the coarser scales.
Fix $d\ge2$, $0<\alpha<1$, and $\beta=1-\alpha$.
We write $B_r(x)$ for the open ball of radius $r$ centered at $x$.
Let $\psi\in C_c^4(\mathbb R^d)$ be nonzero, radial, and supported in
$\overline{B_{1/8}(0)}$. Set
\begin{equation*}
\mathcal A=\{(p,q):1\le p<q\le d\},\qquad
 J^{pq}=e_q\otimes e_p-e_p\otimes e_q,\qquad
 W^\nu=J^\nu\nabla\psi \quad \nu \in \mathcal{A}\, .
\end{equation*}
Thus $W^\nu\in C_c^3$ and $\diver W^\nu=0$.
Following \cite{FM} we define the structure tensor as
\begin{equation*}
 D(z)=\sum_{\nu\in\mathcal A}\int_{\mathbb R^d}
 [W^\nu(y)-W^\nu(y-z)]\otimes[W^\nu(y)-W^\nu(y-z)]\dd y.
\end{equation*}
The tensor is isotropic (cf. \cite{FM}) and for $r>0$ and $|e|=1$ we can decompose it as 
\begin{equation*}
D(re)=D_L(r)e\otimes e+D_N(r)(\Id-e\otimes e).
\end{equation*}
We impose the structural assumption of \cite[Section 1.2]{FM}:
\begin{equation}\label{eq:monotonicity}
 D_L\text{ is nondecreasing on }(0,\infty).
\end{equation}
As pointed out in \cite[Proposition 2.6]{FM} an example of $\psi$ which satisfies all of the above requirements is $\psi(x)=(1-64|x|^2)_+^5$.
But we do not dwell on this here and rather take \eqref{eq:monotonicity} as a standing assumption. 

On the unit torus $\T^d=\mathbb R^d/\mathbb Z^d$, put
\begin{equation}\label{eq:scales}
 \ell_j=2^{-j},\qquad A_j=\ell_j^\alpha,\qquad
 \lambda_j=\ell_j^{-\beta},\qquad
 \tau_j=\kappa 2^{-\lceil\beta j\rceil},\qquad j\ge0.
\end{equation}
The scale-$j$ wavelets are supported in balls of radius $\ell_j/8\le1/8$,
so the construction of \cite[Section 1.2]{FM} starts at scale zero
with these lengths.
The small parameter $0<\kappa\le1$ will be fixed in the proof of the theorem.
In particular,
\begin{equation}\label{eq:lifetimes}
 \frac{\kappa}{2\lambda_j}<\tau_j\le\frac{\kappa}{\lambda_j}.
\end{equation}
Consecutive lifetimes have ratio one or two, and every renewal time
of a coarser scale is a renewal time of every finer scale
\cite[Lemma 2.1]{FM}.

Let $\mathcal I_j=\{0,\ldots,2^j-1\}^d$, $y_i^j=\ell_ji$ modulo
$\mathbb Z^d$, and $Q_j=[0,\ell_j)^d$. Set
\begin{equation}\label{eq:model-horizon}
 T=\kappa\lceil\kappa^{-1}\rceil,\qquad 1\le T\le2.
\end{equation}
On each interval
\begin{equation*}
I_{jk}=[k\tau_j,(k+1)\tau_j),\qquad
 0\le k<T/\tau_j=\lceil\kappa^{-1}\rceil 2^{\lceil\beta j\rceil},
\end{equation*}
choose one shift $\xi_k^j$, uniform on $Q_j$, and coefficients
$\eta_{i,k}^{j,\nu}$, $i\in\mathcal I_j$, $\nu\in\mathcal A$, with
\begin{equation*}
 \E\eta_{i,k}^{j,\nu}=0,\qquad
 \E(\eta_{i,k}^{j,\nu})^2=\sigma^2>0,\qquad
 |\eta_{i,k}^{j,\nu}|\le M.
\end{equation*}
All these individual coefficients and shifts are mutually independent
on one countable product probability space. For every scale $j$, the intervals $I_{jk}$ partition $[0,T)$,
which contains $[0,1)$.

Periodize $W^\nu$ at scale $\ell_j$ by
\begin{equation*}
W_j^\nu(x)=\sum_{n\in\mathbb Z^d}W^\nu((x+n)/\ell_j).
\end{equation*}
The centers and fields are defined inductively in $j$. On $I_{jk}$,
\begin{equation*}
 \dot P_{i,k}^j(t)=u_{<j}(P_{i,k}^j(t),t),\qquad P_{i,k}^j(k\tau_j)=y_i^j+\xi_k^j,
\end{equation*}
\begin{equation*}
 u_j(x,t)=A_j\sum_{i\in\mathcal I_j,\nu\in\mathcal A}\eta_{i,k}^{j,\nu}W_j^\nu(x-P_{i,k}^j(t)),
\end{equation*}
\begin{equation*}
 u_{<j}=\sum_{0\le a<j}u_a.
\end{equation*}
Finally, $u^{(N)}=\sum_{j=0}^Nu_j$ and
$u^{(\infty)}=\sum_{j\ge0}u_j$. By
\cite[Section 1.2 and Lemmas 2.2--2.3]{FM}, the series converges
uniformly for sufficiently small $\kappa$. For some $\kappa_*>0$,
every $0<\kappa\le\kappa_*$ and every realization yield a bounded,
Borel in time, divergence-free field satisfying
\begin{equation}\label{eq:intro-regularity}
 \sup_{0\le t<T}\norm{u^{(\infty)}(\cdot,t)}_{C^\alpha(\T^d)}\le C,
\end{equation}
where $\kappa_*$ and $C$ depend only on $d,\alpha,\psi,M$.
The same estimate holds for every truncated field $u^{(N)}$, uniformly in $N$.
The cited construction estimates are uniform over generations and
apply on $[0,T)$ with the same constants. The bound is
proved in \cite[Lemma 2.2, (2.2)--(2.3)]{FM} (in fact with an argument suggested
by this author, as acknowledged in \cite[Section 1.6, p. 12]{FM}).

\subsection{Main theorem}
For $\mu>0$ and $\theta_0\in L^2(\T^d)$, we let
$\theta^\mu$ be the energy solution of
\begin{equation}\label{eq:PDE}
 \partial_t\theta^\mu+u^{(\infty)}\cdot\nabla\theta^\mu
       =\mu\Delta\theta^\mu,\qquad \theta^\mu(0)=\theta_0.
\end{equation}
Here the solution class is
$L^2(0,T;H^1)\cap H^1(0,T;H^{-1})\subset C([0,T];L^2)$.
We write $\langle f,h\rangle=\int_{\T^d}f\overline h\dd x$ for the
$L^2$ scalar product, and we use the classical Fourier-based norms and seminorms for fractional Sobolev spaces
\begin{equation*}
\norm f_{H^\gamma}^2=\sum_{n\in\mathbb Z^d}
 (1+4\pi^2|n|^2)^\gamma|\widehat f(n)|^2,\qquad
 \norm f_{\dot H^\gamma}^2=\sum_{n\ne0}(2\pi|n|)^{2\gamma}
 |\widehat f(n)|^2.
\end{equation*}
Our main theorem then reads as follows.

\begin{theorem}\label{thm:main}
For every $0<\gamma<\frac{\beta}{2} = \frac{1-\alpha}{2}$ there is $\kappa_\gamma>0$ such that,
for each fixed $0<\kappa\le\kappa_\gamma$ and every deterministic
$\theta_0\in L^2(\T^d)$, the solutions of \eqref{eq:PDE} satisfy
\begin{equation*}
 \sup_{0<\mu\le1}\E\int_0^1\norm{\theta^\mu(t)}_{H^\gamma}^2\dd t
       \le C \norm{\theta_0}_2^2\, ,
\end{equation*}
where the constant depends on $d,\alpha,\psi,M,\sigma,\kappa$, and $\gamma$.
\end{theorem}

\subsection{Acknowledgments and AI declaration} 

The author acknowledges the support of the National Science Foundation
through Grant DMS-2350252 and of the Simons Foundation through the
Simons Initiative on the Geometry of Flows (Grant Award ID
BD-Targeted-00017375).

The author contributed some of the main ideas, but this paper is an outcome
of extensive discussions with ChatGPT (OpenAI), used in its most advanced
commercially available setting; its assistance was fundamental in keeping track of the
most complicated computations and estimates.

\section{Plan of the paper}

From now on we consider the parameter $\gamma$ in Theorem \ref{thm:main} fixed. It is useful to introduce the notation $\theta^{\mu,N}$ for the passive scalar driven
by $u^{(N)}$, namely the solution of 
\begin{equation}\label{eq:PDE-2}
 \partial_t\theta^{\mu,N}+u^{(N)}\cdot\nabla\theta^{\mu,N}
       =\mu\Delta\theta^{\mu, N},\qquad \theta^{\mu,N} (0)=\theta_0.
\end{equation}
To estimate its
$H^\gamma$ seminorm uniformly in $N$ and $\mu$, the main idea is to test the equation with
$g*\theta^{\mu,N}$, with a precise choice of $g$, and take expectations in the random field. For this reason it is convenient to introduce the following quadratic energy.

\begin{definition}[Quadratic energy]
For real even $g\in L^1(\T^d)$ and real $f\in L^2(\T^d)$, set
\begin{equation*}
F_g(f)=\langle g*f,f\rangle
       =\iint f(x)f(\widetilde x)g(x-\widetilde x)\dd x\dd\widetilde x.
\end{equation*}
\end{definition}

$g$ will be chosen so that, in the averaged ``energy identity'' obtained through this testing procedure, a certain portion will yield a coercive estimate of the square of the $L^2_t \dot H^\gamma_x$. This choice of $g$ is the pivotal point of this note and in order to justify it we need to introduce a suitable ``averaged'' version of $D$. 

\begin{definition}[Averaged covariance]
Fix $r_0=1/64$. For $0<|z|<2r_0$, define
\begin{equation}\label{eq:averaged-covariance}
 \mathcal D(z)=\frac{\sigma^2}{2}\sum_{j\ge0}\tau_jA_j^2D(z/\ell_j).
\end{equation}
Write $\mathcal D_L,\mathcal D_N$ for its longitudinal and transverse
components, as for $D$.
\end{definition}

In particular the main point is that, using some of the insights in \cite{FM},
we can choose $g$ in Section~\ref{sec:kernels} so that $\mathcal D:D^2g$
is the fractional Sobolev kernel near the origin. Assumption
\eqref{eq:monotonicity} gives the sign needed for the coercivity. The following is the key lemma which details the choice of $g$.

\begin{lemma}[The adapted radial kernel]\label{lem:g}
There is an even periodic kernel $g\in L^1(\T^d)$, radial on its
support in $\overline{B_{2r_0}(0)}$, $C^4$ off the origin, with
\begin{equation}\label{eq:g-bounds}
 |D^kg(z)|\le C_\gamma |z|^{-d+\beta-2\gamma-k}\quad(0\le k\le4),
\end{equation}
and which satisfies, for $0<r<r_0$, the following pointwise relations
\begin{equation}\label{eq:g-signs}
 \mathcal D:D^2g=\kappa r^{-d-2\gamma},\qquad
 \Delta g>0,\qquad |D^2g|\le C_\gamma\Delta g.
\end{equation}
\end{lemma}

The kernel $g$ is specified in Section~\ref{sec:kernels}. 
As mentioned, the main strategy of the proof is to derive an energy estimate for $F_g (\theta_t^{\mu, N})$ for fixed $\mu$ and $N$, and take its expectation in the random field. Some annoying technicalities must be taken care of, since the kernel $g$ is singular and the identity involves applying differential operators to it. These technicalities are also addressed in Section~\ref{sec:kernels}. 

Section~\ref{sec:renewal} is where the main identity is derived, while in turn Section~\ref{sec:prelim} introduces the tools to justify this derivation. This key identity is
\eqref{eq:renewal} in Proposition \ref{prop:renewal}, but it is, however, more sophisticated than a standard energy identity. Following Fefferman and Mayboroda, cf. \cite{FM}, some terms are rewritten as time integral of further time-derivatives over some suitably selected renewal times. The point is that, thanks to the probabilistic structure, when applying the fundamental theorem of calculus to the relevant quantity, this is found to vanish at the left endpoint of the chosen interval (the expectation and the structure of the model makes this possible, cf. Lemma~\ref{lem:average}). This idea is applied twice and in particular produces a term in the identity \eqref{eq:renewal} (which is denoted by (C)) where \eqref{eq:g-signs} emerges. 

It is then possible to exploit the particular choice of $g$ in Lemma \ref{lem:g} to show that (C) is a coercive term which controls the relevant Sobolev seminorm, modulo some errors (in turn bounded using the classical $L^2$ energy identity). This is the content of Corollary~\ref{c:controllo}: Section~\ref{sec:principal} contains the relevant estimates leading to it. The remaining terms in \eqref{eq:renewal} are treated either as errors, which can be reabsorbed into the Sobolev seminorm choosing the parameter $\kappa$ sufficiently small, or they are controlled using the classical $L^2$ energy identity. The relevant estimates allowing these bounds are collected in Section \ref{sec:estimates}. 

All these ingredients are then combined in Section~\ref{sec:completion} to complete the proof of Theorem \ref{thm:main}. 

\section{Deterministic estimates on the vector fields and generational averages}

We recall here first some simple estimates for the building blocks $u_j$, which are in fact the main point of the regularity bounds on $u^{(N)}$ and $u^{(\infty)}$.
The constants $C$ are independent of $\mu$ and $\gamma$. 

\begin{lemma}\label{lem:spatial}
There is $\kappa_*>0$ such that, for $\kappa\le\kappa_*$, each active
generation has $\ell_j/2$-separated centers, bounded packing in balls
of radius $\ell_j/8$, and
\begin{equation}\label{eq:spatial}
 \norm{D^q u_j}_\infty\le C A_j\ell_j^{-q}\qquad(0\le q\le3).
\end{equation}
\end{lemma}
The spatial assertions are \cite[Lemma 2.2, (2.2)]{FM}.

We next recall another simple fact about the construction, which follows from the relation between the building blocks in the generation $a$ and the building blocks of the previous generations.

\begin{lemma}[A generation average]\label{lem:average}
Condition on any coarser data which determine a measure-preserving
bijection $\Psi:\T^d\to\T^d$. For
$U(x)=A_b\sum_{i,\nu}\eta_i^\nu W_b^\nu(x-\Psi(y_i^b+\xi))$,
let $\E_{\rm new}$ denote integration over the common shift
$\xi$, uniform on $Q_b$, and all coefficients $\eta_i^\nu$ of this
generation, keeping the coarser data and hence $\Psi$ fixed. Then,
when $|z|=|x-\widetilde x|\le1/4$,
\begin{equation*}
 \E_{\rm new}U(x)=0,
\end{equation*}
\begin{equation*}
 \E_{\rm new}(U(x)-U(\widetilde x))^{\otimes2}=\sigma^2A_b^2D(z/\ell_b),
\end{equation*}
\begin{equation*}
 \E_{\rm new}[(U\cdot\nabla)U](x)=0.
\end{equation*}
\end{lemma}
This is \cite[Lemma 6.2]{FM}. Its unfolding condition holds since
$\ell_b/4+|z|\le1/2<1$. The coarser flows used below are
measure-preserving bijections.

\section{Covariance and quadratic kernels}\label{sec:kernels}

Here and below, $C$ depends only on the fixed model, and $C_\gamma$
may also depend on $\gamma$. Unless indicated, these constants are
independent of $\kappa\le\kappa_*$, the cutoffs, $\mu$ and the
localization parameter introduced in Section~\ref{sec:renewal}.
We start by collecting some simple bounds on the tensor $\mathcal{D}$. First of all recall that the tensor $D$ is even and $C^3$, is constant outside $B_{1/4}(0)$, and
satisfies $\diver_zD=0$ and $D_L(r)>0$ for $r>0$; see
\cite[Lemma 2.5 and (2.9)]{FM}. In particular,
\begin{equation}\label{eq:D-bounds}
 |D(z)|\le C\min(|z|^2,1),\qquad
 |D_L'(r)|\le C\min(r,1),\qquad
 |D_L''(r)|+|D_L'''(r)|\le C.
\end{equation}
We also point out that the following relation holds:
\begin{equation}\label{eq:KH}
 D_N(r)=D_L(r)+\frac{rD_L'(r)}{d-1}.
\end{equation}

\begin{lemma}\label{lem:covariance}
For $0<r<2r_0$,
\begin{equation}\label{eq:covariance-KH}
 \mathcal D_N=\mathcal D_L+r\mathcal D_L'/(d-1),
\end{equation}
\begin{equation*}
 \mathcal D_L'\ge0,
\end{equation*}
\begin{equation*}
 c\kappa r^{2-\beta}\le\mathcal D_L(r)\le C\kappa r^{2-\beta},
\end{equation*}
\begin{equation*}
 |\mathcal D_L^{(k)}(r)|\le C\kappa r^{2-\beta-k}\qquad(1\le k\le3).
\end{equation*}
\end{lemma}
\begin{proof}
Let $\ell_n\ge r>\ell_{n+1}$. By \eqref{eq:lifetimes},
$\tau_jA_j^2$ is between fixed positive multiples of
$\kappa\ell_j^{2-\beta}$. Thus
\begin{equation*}
\mathcal D_L(r)\le C\kappa r^2\sum_{j=0}^n\ell_j^{-\beta}
                   +C \kappa \sum_{j>n}\ell_j^{2-\beta}
        \le C\kappa r^{2-\beta}.
\end{equation*}
For the lower bound, $1<r/\ell_{n+1}\le2$ and
$r/2\le\ell_{n+1}<r$, so
\begin{equation*}
 \mathcal D_L(r)\ge\frac{\sigma^2}{2}\tau_{n+1}A_{n+1}^2
 D_L(r/\ell_{n+1})
 \ge c\kappa\ell_{n+1}^{2-\beta}\min_{1\le q\le2}D_L(q)
 \ge c\kappa r^{2-\beta}.
\end{equation*}
Positive-order derivatives of $D_L$
vanish for $r/\ell_j>1/4$, and \eqref{eq:D-bounds} yields
\begin{equation*}
 |\mathcal D_L'(r)|\le C\kappa r\sum_{\ell_j\ge4r}\ell_j^{-\beta}\le C\kappa r^{1-\beta},
\end{equation*}
\begin{equation*}
 |\mathcal D_L^{\prime\prime}(r)|\le C\kappa\sum_{\ell_j\ge4r}\ell_j^{-\beta}\le C\kappa r^{-\beta},
\end{equation*}
\begin{equation*}
 |\mathcal D_L^{\prime\prime\prime}(r)|\le C\kappa\sum_{\ell_j\ge4r}\ell_j^{-1-\beta}\le C\kappa r^{-1-\beta}.
\end{equation*}
These locally finite differentiated sums also justify termwise
differentiation. Finally, \eqref{eq:monotonicity} and \eqref{eq:KH} give
\begin{equation*}
 \mathcal D_L'(r)=\frac{\sigma^2}{2}\sum_j
 \frac{\tau_jA_j^2}{\ell_j}D_L'(r/\ell_j)\ge0,
\end{equation*}
\begin{equation*}
 \mathcal D_N(r)=\frac{\sigma^2}{2}\sum_j\tau_jA_j^2
 \left[D_L(r/\ell_j)+\frac{rD_L'(r/\ell_j)}{(d-1)\ell_j}\right]
 =\mathcal D_L(r)+\frac{r\mathcal D_L'(r)}{d-1}.
\end{equation*}
\end{proof}

\subsection{The adapted kernel}
We are now ready to define the function $g$ which satisfies the requirements of Lemma~\ref{lem:g}. 
The use of a covariance-adapted radial ODE is inspired by
\cite[Section 4.1, Step 1, (4.20)]{DGP}. There the argument concerns a
Gaussian white-in-time field; here the ODE chooses the test kernel.

\begin{definition}[The adapted kernel]\label{def:g}
Consider the unique function $\bar{g}: (0,2r_0) \to \mathbb R$ such that $\bar{g} (r_0) =0$ and  
\begin{equation}\label{eq:g-ode}
 \bar g'(r)=-\frac{\kappa}{2\gamma}\frac{r^{1-d-2\gamma}}{\mathcal D_L(r)}
                  \quad \mbox{for } 0<r<2r_0.
\end{equation}
Then multiply $\bar{g} (|z|)$ by a smooth radial cut-off 
which equals $1$ on $B_{r_0}(0)$ and vanishes outside
$B_{2r_0}(0)$, and periodize to obtain a function on the periodic torus. Denote the resulting difference kernel by $g(x-\widetilde x)$.
By abuse of notation, we use the same letter $g(r)$ for its radial
profile on the coordinate ball, viewed as a function of the real
variable $r>0$ and extended by zero for $r\ge2r_0$.
\end{definition}

\begin{proof}[Proof of Lemma~\ref{lem:g}]
Lemma~\ref{lem:covariance} and differentiation of a reciprocal give
$|(\mathcal D_L^{-1})^{(h)}|\le C\kappa^{-1}r^{\beta-2-h}$ for $0\le h\le3$.
For $1\le q\le4$, differentiation of \eqref{eq:g-ode} gives
\begin{equation}\label{eq:g-radial-derivative-estimate}
 |\bar g^{(q)}(r)| \le C_\gamma\kappa\sum_{m=0}^{q-1} r^{1-d-2\gamma-m}\kappa^{-1} r^{\beta-2-(q-1-m)} \le C_\gamma r^{-d+\beta-2\gamma-q}.
\end{equation}
For $0<r<r_0$, integration from $r$ to $r_0$ yields
\begin{equation*}
 |g(r)|\le C_\gamma\int_r^{r_0}t^{-d+\beta-2\gamma-1}\dd t
 \le C_\gamma r^{-d+\beta-2\gamma}.
\end{equation*}
The cutoff in Definition~\ref{def:g} has fixed derivatives; thus
\eqref{eq:g-radial-derivative-estimate} also bounds $g^{(q)}$.
The Cartesian derivatives satisfy
\begin{equation}\label{eq:g-Cartesian-derivative-estimate}
 |D^qg(z)|\le C_q\sum_{m=1}^qr^{m-q}|g^{(m)}(r)|
 \le C_\gamma r^{-d+\beta-2\gamma-q},\qquad1\le q\le4.
\end{equation}
In particular
\begin{equation*}
\int_{|z|<r_0}|g(z)|\dd z\le C_\gamma
               \int_0^{r_0}r^{\beta-2\gamma-1}\dd r<\infty.
\end{equation*}
By \eqref{eq:covariance-KH}, on $0<r<r_0$,
\begin{equation*}
 \mathcal D:D^2g=\mathcal D_Lg''+(d-1)\mathcal D_Ng'/r =r^{1-d}(r^{d-1}\mathcal D_Lg')' =r^{1-d}\left(-\frac\kappa{2\gamma}r^{-2\gamma}\right)' =\kappa r^{-d-2\gamma},
\end{equation*}
\begin{equation*}
 g''=|g'|\left(\frac{d+2\gamma-1}{r}+\frac{\mathcal D_L'}{\mathcal D_L}\right),
\end{equation*}
\begin{equation*}
 \Delta g=|g'|\left(\frac{2\gamma}{r}+\frac{\mathcal D_L'}{\mathcal D_L}\right).
\end{equation*}
Since $\mathcal D_L'\ge0$, the tangential eigenvalue satisfies
$|g'|/r\le\Delta g/(2\gamma)$, while
$g''=\Delta g+(d-1)|g'|/r$. Consequently,
\begin{equation*}
 |D^2g|\le g''+(d-1)|g'|/r
 \le C_\gamma\Delta g,
\end{equation*}
which completes \eqref{eq:g-signs}.
Finally $\mathcal D_L/\kappa$ is independent of $\kappa$ by \eqref{eq:scales}.
\end{proof}

\subsection{Linear and bilinear operators}
Now that we have defined the integrable energy kernel $g$, we wish
to consider kernels obtained by applying spatial differential
operators to it. These kernels are less regular and may also depend
on the centers in the construction, which we treat as independent
variables. We therefore specify their linear and bilinear actions
before evaluating the centers along their paths.

\begin{definition}[Configuration kernels]\label{def:kernel}
Fix finitely many active generations. The variables
$(x,\widetilde x;P)$ consist of two independent torus coordinates and
independent center coordinates $P=(P_i^m)$. A real symmetric
configuration kernel satisfies
$H(x,\widetilde x;P)=H(\widetilde x,x;P)$. Bounded absolute row and
column integrals, uniformly over the center configurations under
consideration, mean
\begin{equation*}
 \norm H_{\rm Schur}:=\max\left\{
 \sup_P\mathop{\rm ess\,sup}_{x\in\T^d}
       \int_{\T^d}|H(x,\widetilde x;P)|\dd\widetilde x,\ 
 \sup_P\mathop{\rm ess\,sup}_{\widetilde x\in\T^d}
       \int_{\T^d}|H(x,\widetilde x;P)|\dd x\right\}<\infty.
\end{equation*}
For such a kernel the linear operator on $L^2(\T^d)$ is
\begin{equation*}
 Hf(x)=\int_{\T^d}H(x,\widetilde x;P)f(\widetilde x)\dd\widetilde x.
\end{equation*}
On real-valued functions its bilinear operator and quadratic specialization are
\begin{equation*}
 F_H(f,h)=\langle Hf,h\rangle,\qquad F_H(f)=F_H(f,f).
\end{equation*}
Schur's test gives
\begin{equation*}
\norm{Hf}_2\le\norm H_{\rm Schur}\norm f_2;
\end{equation*}
see \cite[Theorems 8.3.1 and 8.3.4(b), pp. 363--365]{Heil}, applied
on a unit fundamental cube. For distributional kernels we specify separately the domains of
their linear and bilinear operators. A radial kernel
supported in $\overline{B_{2r_0}(0)}$ means the periodic difference kernel
$g(x-\widetilde x)$, with the shortest representative used on its support.
\end{definition}

\begin{definition}[Spatial operators]
Following \cite[Definition 4.2]{FM}, fix integers $N\ge J\ge0$
and let $v=u^{(N)}-u^{(J)}$. For $0\le a\le J$, define on
independent coordinates
\begin{equation}\label{eq:B}
 \mathcal B_a=u_a(x)\cdot\nabla_x+u_a(\widetilde x)\cdot\nabla_{\widetilde x}+\sum_{a<m\le J,i}u_a(P_i^m)\cdot\nabla_{P_i^m},
\end{equation}
\begin{equation*}
 \mathcal B_v=v(x)\cdot\nabla_x+v(\widetilde x)\cdot\nabla_{\widetilde x},
\end{equation*}
\begin{equation*}
 \Delta_Y=\Delta_x+\Delta_{\widetilde x}.
\end{equation*}
Here $u_a(y)$ is the active wavelet expression with its coefficients
fixed and its centers independent. The operator $\mathcal B_v$
acts only on the two particle coordinates, with no center derivative.
Operators are applied before the centers are evaluated at their
current positions. The sum in \eqref{eq:B} differentiates the centers of generation
$m$ along their velocity $u_a$, for $a<m$.
If $H=H(x,\widetilde x;P,y)$ also depends on a separate spatial
variable $y$, the convention is
\begin{equation*}
 \mathcal B_aH=u_a(x)\cdot\nabla_xH
 +u_a(\widetilde x)\cdot\nabla_{\widetilde x}H
 +\sum_{a<m\le J,i}u_a(P_i^m)\cdot\nabla_{P_i^m}H
 +u_a(y)\cdot\nabla_yH.
\end{equation*}
In particular, $u_b(y;P)$ depends on $y$ and on the scale-$b$
centers, but not on $x,\widetilde x$. Differentiating these variables
independently gives
\begin{equation}\label{eq:B-ub}
 \mathcal B_au_b(y)
 =A_b\sum_{i,\nu}\eta_i^{b,\nu}
 \bigl[u_a(y)-\ind_{\{a<b\}}u_a(P_i^b)\bigr]\cdot
                                  \nabla W_b^\nu(y-P_i^b).
\end{equation}
The factor $-\ind_{\{a<b\}}u_a(P_i^b)$ in \eqref{eq:B-ub}
accounts for the motion of $P_i^b$ under $u_a$ when $a<b$. Values at $y=x$ or
$y=\widetilde x$ are taken only after these derivatives have been
computed. This is the meaning of $\mathcal B_au_b(x)$ and
$\mathcal B_au_b(\widetilde x)$ below.
\end{definition}

We are now ready to define the kernels induced by these operators acting on $g$.

\begin{definition}[Regularized kernels]\label{def:responses}
Choose a smooth scalar cutoff $\eta$ equal to zero on $[0,1]$ and one
on $[2,\infty)$. For $0<2\varepsilon<r_0$, set
\begin{equation*}
g_\varepsilon(r)=\eta(r/\varepsilon)g(r)
                      +(1-\eta(r/\varepsilon))g(\varepsilon)
\end{equation*}
for $0<r<2r_0$. Set $g_\varepsilon(0)=g(\varepsilon)$ and
$g_\varepsilon=g$ outside $B_{2r_0}(0)$. These are periodic
$C^4$ kernels, constant near zero. Define $\mathcal B_ag$ and, at
fixed $N$, $\mathcal B_vg$ as the Schur-norm limits of
$\mathcal B_ag_\varepsilon$ and $\mathcal B_vg_\varepsilon$.
Define $\Delta_Yg$ as the distributional limit of
$\Delta_Yg_\varepsilon$. As a difference distribution it is
$2\Delta_zg$, and its action on $L^2$ is
$(\Delta_Yg)*f=2\Delta(g*f)$ in distributions. Its linear
action $H^1\to H^{-1}$ and its bilinear action on real $H^1$ functions are
\begin{equation}\label{eq:lap-action}
 F_{\Delta_Yg}(f,h)
 :=\langle(\Delta_Yg)*f,h\rangle_{H^{-1},H^1}
 =-2\langle\nabla(g*f),\nabla h\rangle,\qquad f,h\in H^1.
\end{equation}
Its quadratic specialization is
\begin{equation*}
 F_{\Delta_Yg}(f)=-2\langle\nabla f,\nabla(g*f)\rangle.
\end{equation*}
\end{definition}

The distribution $\Delta_Yg$ is not an $L^1$ kernel and does not
define a bounded operator on $L^2$. Its linear and bilinear actions
therefore require the domains in Definition~\ref{def:responses}.
The following lemma justifies that definition.

\begin{lemma}[Convergence and cancellation]\label{lem:pairings}
Let $g_\varepsilon$ be defined using $\eta$ as in
Definition~\ref{def:responses}. The following assertions hold;
constants in estimates involving $g_\varepsilon$ may depend on
this fixed auxiliary function.
\begin{enumerate}
\renewcommand{\theenumi}{\roman{enumi}}
\renewcommand{\labelenumi}{(\theenumi)}
\item\label{item:g-convergence} The kernels $g_\varepsilon$ converge to
$g$ in $L^1$, and their associated linear and bilinear operators
converge on $L^2$. More precisely,
\begin{equation*}
 \norm{g_\varepsilon-g}_1\le C_\gamma\varepsilon^{\beta-2\gamma}.
\end{equation*}
For real $f,h\in L^2$,
\begin{equation*}
 |F_g(f,h)|\le\norm g_1\norm f_2\norm h_2.
\end{equation*}
\item\label{item:B-convergence} At each fixed scale the Schur limit is
\begin{equation}\label{eq:B-g-increment}
 \mathcal B_ag(x,\widetilde x)
 =[u_a(x)-u_a(\widetilde x)]\cdot\nabla g(x-\widetilde x).
\end{equation}
It satisfies
\begin{equation*}
 \norm{\mathcal B_ag}_{\rm Schur}\le C_\gamma\lambda_a,
\end{equation*}
\begin{equation}\label{eq:B-g-error}
 \norm{\mathcal B_ag_\varepsilon-\mathcal B_ag}_{\rm Schur}
 \le C_\gamma\lambda_a\varepsilon^{\beta-2\gamma}.
\end{equation}
The same assertions hold for $\mathcal B_v$, with $u_a$ replaced
by $v$ and $\lambda_a$ by $\norm{Dv}_\infty$, at fixed $N$.
These kernels define bounded linear and bilinear operators on $L^2$.
\item\label{item:lap-convergence} The operators $\Delta_Yg_\varepsilon$
converge to $2\Delta(g*\cdot)$ in operator norm from $L^2$ to
$H^{-2}$ and from $H^1$ to $H^{-1}$. The limiting bilinear
operator is \eqref{eq:lap-action}, and
\begin{equation*}
 |F_{\Delta_Yg}(f,h)|\le2\norm g_1\norm{\nabla f}_2\norm{\nabla h}_2,
 \qquad f,h\in H^1.
\end{equation*}
\item\label{item:cancellation} The operators
$\mathcal B_a,\mathcal B_v,\Delta_Y$ preserve zero row and column
integrals of regular configuration kernels. Each sends
$g_\varepsilon$ to a kernel with zero row and column integrals.
For a symmetric kernel $H$ with these cancellations and finite
Schur norm,
\begin{equation}\label{eq:increments}
 F_H(f)=-\frac12\iint|f(x)-f(\widetilde x)|^2H(x,\widetilde x)
                  \dd x\dd\widetilde x,\qquad f\in L^2(\T^d;\mathbb R).
\end{equation}
\item\label{item:lap-increments} For real $h\in H^1$,
\begin{equation}\label{eq:pairing-bounds}
 F_{\Delta_Yg}(h)=-\int_{\T^d}\norm{h(\cdot+z)-h}_2^2\Delta_zg(z)\dd z.
\end{equation}
The integral uses the off-diagonal function $\Delta_zg$ and is
absolutely convergent. Moreover,
\begin{equation}\label{eq:lap-increment-limit}
 F_{\Delta_Yg}(h)=\lim_{\varepsilon\downarrow0}F_{\Delta_Yg_\varepsilon}(h) =-\frac12\lim_{\varepsilon\downarrow0}\iint |h(x)-h(\widetilde x)|^2(\Delta_x+\Delta_{\widetilde x}) g_\varepsilon(x-\widetilde x)\dd x\dd\widetilde x.
\end{equation}
These limits are uniform on bounded subsets of $H^1$.
\item\label{item:eta-independence} All limiting linear and bilinear
operators in (i)--(v) are independent of the auxiliary function
$\eta$ in Definition~\ref{def:responses}.
\end{enumerate}
\end{lemma}
\begin{proof}
We prove (\ref{item:g-convergence}). The support of
$g_\varepsilon-g$ and \eqref{eq:g-bounds} give
\begin{equation*}
 \norm{g_\varepsilon-g}_1 \le C_\gamma\int_0^{2\varepsilon}r^{\beta-2\gamma-1}\dd r +C_\gamma\varepsilon^{-d+\beta-2\gamma}(2\varepsilon)^d \le C_\gamma\varepsilon^{\beta-2\gamma}.
\end{equation*}
Young's inequality and Cauchy--Schwarz give, respectively,
\begin{equation*}
 \norm{(g_\varepsilon-g)*f}_2
 \le\norm{g_\varepsilon-g}_1\norm f_2
 \le C_\gamma\varepsilon^{\beta-2\gamma}\norm f_2,
\end{equation*}
\begin{equation*}
 |F_{g_\varepsilon}(f,h)-F_g(f,h)|
 \le\norm{g_\varepsilon-g}_1\norm f_2\norm h_2
 \le C_\gamma\varepsilon^{\beta-2\gamma}\norm f_2\norm h_2.
\end{equation*}
The same inequalities yield
\begin{equation*}
 |F_g(f,h)|\le\norm{g*f}_2\norm h_2
                    \le\norm g_1\norm f_2\norm h_2.
\end{equation*}

We prove (\ref{item:B-convergence}). The increment $|u_a(x)-u_a(\widetilde x)|\le C\lambda_a r$
makes the kernel in \eqref{eq:B-g-increment} integrable. Its inner row integral and
that of its regularization are bounded by
\begin{equation*}
C_\gamma\lambda_a\int_0^{2\varepsilon}
                   r^{\beta-2\gamma-1}\dd r
       \le C_\gamma\lambda_a\varepsilon^{\beta-2\gamma}
\end{equation*}
Symmetry gives the column estimate. Integration over the full support
also gives
\begin{equation*}
 \norm{\mathcal B_ag}_{\rm Schur}
 \le C_\gamma\lambda_a\int_0^{2r_0}r^{\beta-2\gamma-1}\dd r
 \le C_\gamma\lambda_a.
\end{equation*}
Schur's test \cite[Theorems 8.3.1 and 8.3.4(b)]{Heil}, on a
unit fundamental cube, now gives
\begin{equation*}
 \norm{(\mathcal B_ag_\varepsilon-\mathcal B_ag)f}_2
 \le C_\gamma\lambda_a\varepsilon^{\beta-2\gamma}\norm f_2,
\end{equation*}
\begin{equation*}
 |F_{\mathcal B_ag_\varepsilon}(f,h)-F_{\mathcal B_ag}(f,h)|
 \le C_\gamma\lambda_a\varepsilon^{\beta-2\gamma}\norm f_2\norm h_2.
\end{equation*} For $\mathcal B_vg$ use
$|v(x)-v(\widetilde x)|\le\norm{Dv}_\infty|x-\widetilde x|$
in the same integral, which gives the asserted Schur limits at fixed $N$.

We prove (\ref{item:lap-convergence}). Young's inequality gives
\begin{equation*}
 \norm{2\Delta[(g_\varepsilon-g)*f]}_{H^{-2}}
 \le2\norm{g_\varepsilon-g}_1\norm f_2,\qquad f\in L^2,
\end{equation*}
\begin{equation*}
 \norm{2\Delta[(g_\varepsilon-g)*f]}_{H^{-1}}
 \le2\norm{g_\varepsilon-g}_1\norm{\nabla f}_2,\qquad f\in H^1.
\end{equation*}
Together with (\ref{item:g-convergence}), these prove the two
operator limits and \eqref{eq:lap-action}. Moreover,
\begin{equation*}
 |F_{\Delta_Yg}(f,h)|
 =2|\langle\nabla(g*f),\nabla h\rangle|
 \le2\norm g_1\norm{\nabla f}_2\norm{\nabla h}_2.
\end{equation*}

We prove (\ref{item:cancellation}). For a regular $H$ put
$R_H(x;P)=\int H(x,\widetilde x;P)\dd\widetilde x$. Then
\begin{equation}\label{eq:row-cancellation}
\int\mathcal B_aH\dd\widetilde x
 =u_a(x)\cdot\nabla_xR_H
  +\sum_{a<m\le J,i}u_a(P_i^m)\cdot\nabla_{P_i^m}R_H
  +\int\diver_{\widetilde x}(u_a(\widetilde x)H)\dd\widetilde x.
\end{equation}
The integral of $\diver_{\widetilde x}(u_a(\widetilde x)H)$ in
\eqref{eq:row-cancellation} vanishes by periodicity. Thus zero row integrals are preserved,
and a constant center-independent row integral is sent to zero.
The same computation gives
\begin{equation*}
 \int\mathcal B_vH\dd\widetilde x=v(x)\cdot\nabla_xR_H,
\end{equation*}
\begin{equation*}
 \int\Delta_YH\dd\widetilde x=\Delta_xR_H.
\end{equation*}
Symmetry gives the assertions for column integrals. If these
integrals vanish, expansion of the square yields
\begin{equation*}
 \iint|f(x)-f(\widetilde x)|^2H\dd x\dd\widetilde x
 =\int f(x)^2R_H(x)\dd x
   +\int f(\widetilde x)^2\left(\int H\dd x\right)\dd\widetilde x
   -2F_H(f) =-2F_H(f).
\end{equation*}
The finite Schur norm makes every integral absolutely convergent,
so this proves \eqref{eq:increments}.

We prove (\ref{item:lap-increments}), including the asserted
uniformity under $H^1$ approximation.
In the rest of this proof write $s=\beta-2\gamma\in(0,1)$ and take
$0<2\varepsilon<r_0$. Start with a smooth real $h$. The fundamental
theorem of calculus and Cauchy--Schwarz give, for $r=|z|$,
\begin{equation*}
 h(x+z)-h(x)=\int_0^1z\cdot\nabla h(x+tz)\dd t,
\end{equation*}
\begin{equation}\label{eq:H1-translation}
 \norm{h(\cdot+z)-h}_2^2\le r^2\int_0^1\int_{\T^d}|\nabla h(x+tz)|^2\dd x\dd t=r^2\norm{\nabla h}_2^2.
\end{equation}
On $\varepsilon<r<2\varepsilon$, put $\eta_\varepsilon(r)=\eta(r/\varepsilon)$.
The radial Laplacian of the function in Definition~\ref{def:responses} is
\begin{equation}\label{eq:regularized-radial-laplacian}
 \Delta g_\varepsilon =\eta_\varepsilon\Delta g +\frac2\varepsilon\eta'(r/\varepsilon)g' +\left[\frac1{\varepsilon^2}\eta''(r/\varepsilon) +\frac{d-1}{\varepsilon r}\eta'(r/\varepsilon)\right] [g(r)-g(\varepsilon)].
\end{equation}
By \eqref{eq:g-bounds},
\begin{equation}\label{eq:g-cutoff-difference}
|g(r)-g(\varepsilon)|\le\int_\varepsilon^{2\varepsilon}|g'(t)|\dd t
 \le C_\gamma\int_\varepsilon^{2\varepsilon}t^{-d+s-1}\dd t
 \le C_\gamma\varepsilon^{-d+s}.
\end{equation}
On this annulus, \eqref{eq:g-bounds} and
\eqref{eq:g-cutoff-difference} give
\begin{equation*}
 |\Delta g_\varepsilon|
 \le|\Delta g|+C_\eta\varepsilon^{-1}|g'|
       +C_\eta\varepsilon^{-2}|g(r)-g(\varepsilon)|
 \le C_\gamma\varepsilon^{-d+s-2}.
\end{equation*} For $r\le\varepsilon$ the
Laplacian is zero, and for $r\ge2\varepsilon$ it equals $\Delta g$.
Also $|\Delta g(z)|\le C_\gamma r^{-d+s-2}$ on $0<r<r_0$.
Using polar coordinates and \eqref{eq:H1-translation}, we obtain
separately
\begin{equation}\label{eq:inner-laplacians}
 \int_{|z|<2\varepsilon}\norm{h(\cdot+z)-h}_2^2|\Delta g_\varepsilon(z)|\dd z \le C_\gamma\norm{\nabla h}_2^2
 \varepsilon^{-d+s-2}\int_0^{2\varepsilon}r^{d+1}\dd r
 \le C_\gamma\varepsilon^s\norm{\nabla h}_2^2,
\end{equation}
\begin{equation}\label{eq:inner-unregularized-laplacian}
 \int_{|z|<2\varepsilon}\norm{h(\cdot+z)-h}_2^2|\Delta g(z)|\dd z \le C_\gamma\norm{\nabla h}_2^2\int_0^{2\varepsilon}r^{s-1}\dd r \le C_\gamma\varepsilon^s\norm{\nabla h}_2^2.
\end{equation}
The same calculation up to $r_0$, and boundedness on the fixed
outer annulus, show absolute convergence of the off-diagonal
integral in \eqref{eq:pairing-bounds}. Since the kernels agree
outside $B_{2\varepsilon}(0)$, \eqref{eq:inner-laplacians} and \eqref{eq:inner-unregularized-laplacian} give
\begin{equation}\label{eq:lap-form-error}
 \left|\int_{\T^d}\norm{h(\cdot+z)-h}_2^2
                       [\Delta g_\varepsilon(z)-\Delta g(z)]\dd z\right|
       \le C_\gamma\varepsilon^s\norm{\nabla h}_2^2.
\end{equation}
Here $\Delta g$ in the integral is its off-diagonal function.

For the regularized kernel,
$(\Delta_x+\Delta_{\widetilde x})g_\varepsilon(x-\widetilde x)
=2\Delta_zg_\varepsilon(x-\widetilde x)$. Periodic integration by
parts therefore yields
\begin{equation*}
 F_{\Delta_Yg_\varepsilon}(h)
 =2\int_{\T^d}h(x)\Delta(g_\varepsilon*h)(x)\dd x =-2\sum_{i=1}^d\int_{\T^d}\partial_i h(x)
                          \partial_i(g_\varepsilon*h)(x)\dd x
 =-2\langle\nabla h,\nabla(g_\varepsilon*h)\rangle.
\end{equation*}
Moreover $\int_{\T^d}\Delta g_\varepsilon=0$. Expanding the square
and changing variables gives
\begin{equation*}
 \begin{aligned}
 &\int_{\T^d}\norm{h(\cdot+z)-h}_2^2\Delta g_\varepsilon(z)\dd z =\iint[h(x)^2+h(\widetilde x)^2-2h(x)h(\widetilde x)]
              \Delta g_\varepsilon(x-\widetilde x)\dd x\dd\widetilde x\\
 &=2\norm h_2^2\int_{\T^d}\Delta g_\varepsilon(z)\dd z
                    -2\langle h,\Delta(g_\varepsilon*h)\rangle
 =-F_{\Delta_Yg_\varepsilon}(h).
 \end{aligned}
\end{equation*}
On the other hand, Young's inequality
and the $L^1$ convergence already proved give
\begin{equation*}
 \left|F_{\Delta_Yg_\varepsilon}(h)
              +2\langle\nabla h,\nabla(g*h)\rangle\right| =2\left|\langle\nabla h,(g_\varepsilon-g)*\nabla h\rangle\right|
 \le2\norm{g_\varepsilon-g}_1\norm{\nabla h}_2^2
 \le C_\gamma\varepsilon^s\norm{\nabla h}_2^2.
\end{equation*}
Together with \eqref{eq:lap-form-error}, this proves the Laplacian
identities in \eqref{eq:pairing-bounds} and
\eqref{eq:lap-increment-limit} for smooth $h$.

These passages are uniform under $H^1$ approximation. If $h,k$ are smooth, Cauchy--Schwarz and
\eqref{eq:H1-translation} give
\begin{equation*}
 \left|\norm{h(\cdot+z)-h}_2^2-\norm{k(\cdot+z)-k}_2^2\right| \le r^2\norm{\nabla(h-k)}_2 (\norm{\nabla h}_2+\norm{\nabla k}_2).
\end{equation*}
The weighted integral needed here satisfies
\begin{equation}\label{eq:lap-weighted-moment}
 \int_{\T^d}|z|^2|\Delta g(z)|\dd z\le C_\gamma\int_0^{r_0}r^{s-1}\dd r+C_\gamma\le C_\gamma,
\end{equation}
\begin{equation}\label{eq:lap-regularized-moment}
 \sup_{0<2\varepsilon<r_0}\int_{\T^d}|z|^2|\Delta g_\varepsilon(z)|\dd z \le C_\gamma.
\end{equation}
Identify $\T^d$ with $[-1/2,1/2)^d$ in these integrals. The
off-diagonal functions $\Delta g$ and $\Delta g_\varepsilon$ are
supported in $\overline{B_{2r_0}(0)}$, with $2r_0=1/32<1/2$.
Every torus point in that ball has a unique representative in the
Euclidean ball, and $|z|$ denotes its Euclidean norm. Both integrands
are zero outside the ball, so their weighted moments are precisely
the Euclidean radial integrals in \eqref{eq:lap-weighted-moment}
and \eqref{eq:lap-regularized-moment}. For either off-diagonal function $K=\Delta g$ or
$K=\Delta g_\varepsilon$, the translation bound \eqref{eq:H1-translation}
and these weighted moments give
\begin{equation*}
 \left|\int_{\T^d}\left[\norm{h(\cdot+z)-h}_2^2
 -\norm{k(\cdot+z)-k}_2^2\right]K(z)\dd z\right|
 \le C_\gamma\norm{\nabla(h-k)}_2(\norm{\nabla h}_2+\norm{\nabla k}_2),
\end{equation*}
uniformly in $\varepsilon$.
The integration-by-parts forms obey
\begin{equation*}
 \left|\langle\nabla h,\nabla(g_\varepsilon*h)\rangle -\langle\nabla k,\nabla(g_\varepsilon*k)\rangle\right| \le\norm{g_\varepsilon}_1\norm{\nabla(h-k)}_2 (\norm{\nabla h}_2+\norm{\nabla k}_2),
\end{equation*}
and $\sup_\varepsilon\norm{g_\varepsilon}_1<\infty$.
Taking a smooth approximating sequence $h_m\to h$ in $H^1$ thus extends every identity and
\eqref{eq:lap-form-error} to $h\in H^1$. The error is uniform on
bounded subsets of $H^1$.

We prove (\ref{item:eta-independence}). Let $\widetilde\eta$
be another auxiliary function as in Definition~\ref{def:responses},
and denote by $\widetilde g_\delta$ its regularization. The
calculation in \eqref{eq:regularized-radial-laplacian} applies with
$\widetilde\eta,\delta$ in place of $\eta,\varepsilon$ and gives
\eqref{eq:inner-laplacians} with a constant depending on
$\widetilde\eta$. Both regularizations converge to $g$ in $L^1$,
so they give the same limiting linear and bilinear operators.
In particular,
\begin{equation*}
 |F_{\Delta_Yg_\varepsilon}(h)
       -F_{\Delta_Y\widetilde g_\delta}(h)|
 \le C_{\gamma,\eta,\widetilde\eta}
       (\varepsilon^{\beta-2\gamma}+\delta^{\beta-2\gamma})\norm{\nabla h}_2^2.
\end{equation*}
Likewise, \eqref{eq:B-g-error} for the two auxiliary functions gives
\begin{equation*}
 \norm{\mathcal B_ag_\varepsilon-\mathcal B_a\widetilde g_\delta}_{\rm Schur}
 \le C_{\gamma,\eta,\widetilde\eta}\lambda_a
       (\varepsilon^{\beta-2\gamma}+\delta^{\beta-2\gamma}),
\end{equation*}
and the same estimate holds for $\mathcal B_v$ with
$\lambda_a$ replaced by $\norm{Dv}_\infty$.
\end{proof}

\subsection{Classical estimates on Sobolev spaces}
We shall use the following classical relations between the Fourier norms and seminorms on fractional Sobolev spaces and the Gagliardo-type norms, which we recall for the reader's convenience:
\begin{equation}\label{eq:fractional-translation}
 \norm{f(\cdot+z)-f}_2\le C_\gamma|z|^\gamma\norm f_{\dot H^\gamma},
\end{equation}
\begin{equation}\label{eq:fractional-upper}
 \int_{|z|<2r_0}\frac{\norm{f(\cdot+z)-f}_2^2}{|z|^{d+2\gamma}}\dd z\le C_\gamma\norm f_{\dot H^\gamma}^2,
\end{equation}
\begin{equation}\label{eq:fractional}
 \norm f_{\dot H^\gamma}^2\le C_\gamma\left(\int_{|z|<r_0}\frac{\norm{f(\cdot+z)-f}_2^2}{|z|^{d+2\gamma}}\dd z+\norm f_2^2\right).
\end{equation}
Parseval's identity gives the computation
\begin{equation*}
 \norm{f(\cdot+z)-f}_2^2 =\sum_n|e^{2\pi in\cdot z}-1|^2|\widehat f(n)|^2 \le C_\gamma|z|^{2\gamma} \sum_{n\ne0}(2\pi|n|)^{2\gamma}|\widehat f(n)|^2 =C_\gamma|z|^{2\gamma}\norm f_{\dot H^\gamma}^2.
\end{equation*}
For the other two, \cite[Theorem 1.5, (1.6)--(1.7)]{RS}, rescaled
from the $2\pi$ torus and symmetrized, gives
\begin{equation*}
\norm f_{\dot H^\gamma}^2
 =\frac{4^\gamma\Gamma(d/2+\gamma)}{2\pi^{d/2}|\Gamma(-\gamma)|}
 \int_{\T^d}\norm{f(\cdot+z)-f}_2^2
                  \sum_{n\in\mathbb Z^d}|z+n|^{-d-2\gamma}\dd z.
\end{equation*}
The periodic kernel is comparable with $|z|^{-d-2\gamma}$ near zero
and bounded on its complement. Fourier approximation and Tonelli's
theorem extend the identity to $H^\gamma$.

\section{Scalar and quadratic energy identities}\label{sec:prelim}

We use the scalar energy identity \cite[Lemma 8.2(a)]{FM} and
the two-particle evolution supplied by \cite[Lemmas 8.6(a) and
8.2(c), pp. 56--57]{FM}. The center chain rule then gives the
configuration-kernel identity needed on each renewal interval.

Throughout the paper, $T$ is the endpoint fixed in \eqref{eq:model-horizon}.
The identities in this section also hold on every subinterval $[0,S]\subset[0,T]$. 
For a bounded Borel divergence-free drift, existence in the energy
class, the continuous $L^2$ trace and uniqueness are supplied by
\cite[Section 8, (E1)--(E2) and Lemma 8.2(a)]{FM}. These hypotheses
hold for every finite field and for $u^{(\infty)}$. The same lemma
gives
\begin{equation}\label{eq:energy}
 \norm{\theta_S}_2^2+2\mu\int_0^S\norm{\nabla\theta_t}_2^2\dd t
       =\norm{\theta_0}_2^2\qquad(0\le S\le T),
\end{equation}
and $L^2$ contraction for solutions with the same drift.
Measurability for the truncated fields $u^{(N)}$ is \cite[Lemma 9.1]{FM}, whose
last assertion applies to every positive $\mu$; the measurability for the fixed-viscosity
limit in the case of $u^{(\infty)}$ is justified in Section~\ref{sec:completion}.

\begin{lemma}\label{lem:product}
Fix integers $N\ge J\ge0$ and $\mu>0$, let
$\theta=\theta^{\mu,N}$ solve \eqref{eq:PDE-2} with
$\theta_0\in L^2(\T^d)$, and put
$v=u^{(N)}-u^{(J)}$. Let $K$ be a symmetric $C^2$ configuration
kernel involving only scales $0,\ldots,J$. An active generation
participates in $K$ if its coefficients or shift, or any of its
center coordinates $P_i^m$, occur in $K$. Assume that $K$ has no
explicit time dependence while these generation data are fixed.
Between renewals of participating generations, $F_K(\theta_t)$ is
absolutely continuous and satisfies almost everywhere
\begin{equation}\label{eq:product}
 \frac\dd{\dd t}F_K(\theta_t)
   =\sum_{a=0}^JF_{\mathcal B_aK}(\theta_t)
       +F_{\mathcal B_vK}(\theta_t)
       +\mu F_{\Delta_YK}(\theta_t).
\end{equation}
\end{lemma}

For the reader's convenience, we write \eqref{eq:product} in integral
form. With $P(t)$ denoting the active centers, it reads
\begin{equation}\label{eq:product-integral}
 \begin{aligned}
 &\frac\dd{\dd t}\iint_{\T^d\times\T^d}
       \theta_t(x)\theta_t(\widetilde x)
                     K(x,\widetilde x;P(t))\dd x\dd\widetilde x\\
 &=\iint_{\T^d\times\T^d}\theta_t(x)\theta_t(\widetilde x)
 \Bigl[u^{(N)}(x,t)\cdot\nabla_xK
       +u^{(N)}(\widetilde x,t)\cdot\nabla_{\widetilde x}K\\
 &\hspace{38mm}+\mu(\Delta_x+\Delta_{\widetilde x})K
       +\sum_{m=0}^J\sum_{i\in\mathcal I_m}
             u_{<m}(P_i^m,t)\cdot\nabla_{P_i^m}K\Bigr]_{P=P(t)}
                                               \dd x\dd\widetilde x.
 \end{aligned}
\end{equation}
The derivatives of $K$ are taken before evaluating the centers. A comparison between \eqref{eq:product-integral} with \eqref{eq:product} should justify to the reader the reason for introducing the shorthand notation $F_H (\theta_t)$. In fact we will handle identities which fully expanded would be much longer than \eqref{eq:product-integral}. 

\begin{proof}
By \cite[Lemma 8.6(a)]{FM}, $\theta\otimes\theta$ is the energy
solution of the two-particle equation. Subdivide at the finite-field
renewals: the center paths are then $C^1$, so a $C^2$ configuration
kernel has the time and particle regularity required by
\cite[Lemma 8.2(c)]{FM}. Apply that identity to $K$ evaluated at
the moving centers. This gives the particle terms in
\eqref{eq:product-integral}. The additional time derivative is the
center chain rule,
\begin{equation*}
\partial_tK=\sum_{m\le J,i}u_{<m}(P_i^m)\cdot\nabla_{P_i^m}K
 =\sum_{a\le J}\sum_{a<m\le J,i}
                     u_a(P_i^m)\cdot\nabla_{P_i^m}K.
\end{equation*}
Together with $u^{(N)}=\sum_{a\le J}u_a+v$, this gives
\eqref{eq:product}. The operator $\mathcal B_v$ has no center derivative.
The scalar trace is continuous at every renewal.
\end{proof}

\begin{corollary}\label{cor:singular-product}
Fix integers $N\ge J\ge0$, $\mu>0$, and a datum
$\theta_0\in L^2(\T^d)$. Let $\theta=\theta^{\mu,N}$ solve
\eqref{eq:PDE-2} and put $v=u^{(N)}-u^{(J)}$. Then, for every
$0\le S\le T$,
\begin{equation}\label{eq:singular-energy}
 F_g(\theta_S)-F_g(\theta_0)
 =\int_0^S\left\{\sum_{a=0}^JF_{\mathcal B_ag}(\theta_t)
                   +F_{\mathcal B_vg}(\theta_t)
                   +\mu F_{\Delta_Yg}(\theta_t)\right\}\dd t.
\end{equation}
The bilinear operators associated with $\mathcal B_ag$ and
$\mathcal B_vg$ use the Schur kernels of
Definition~\ref{def:responses}; the Laplacian form uses the
$H^1$ pairing in \eqref{eq:lap-action}. All terms are integrable.
\end{corollary}
\begin{proof}
The kernel $g_\varepsilon$ has no coefficients, shifts or center
variables. Apply
Lemma~\ref{lem:product} on the intervals between the finitely many
renewals of $u^{(N)}$. At a renewal time $t$, the two kernel traces
are the same and $\theta$ has a continuous $L^2$ trace. More
explicitly, for $r,t\in[0,T]$,
\begin{equation*}
|F_{g_\varepsilon}(\theta_r)-F_{g_\varepsilon}(\theta_t)|
 \le\norm{g_\varepsilon}_1\norm{\theta_r-\theta_t}_2
                           (\norm{\theta_r}_2+\norm{\theta_t}_2),
\end{equation*}
so the endpoint terms match as $r\to t$. The intervalwise
absolutely continuous forms therefore concatenate to an absolutely
continuous function on $[0,T]$. Integrating its derivative gives
\begin{equation}\label{eq:regularized-energy}
 F_{g_\varepsilon}(\theta_S)-F_{g_\varepsilon}(\theta_0)
 =\int_0^S\left\{\sum_{a=0}^JF_{\mathcal B_ag_\varepsilon}(\theta_t)
              +F_{\mathcal B_vg_\varepsilon}(\theta_t)
              +\mu F_{\Delta_Yg_\varepsilon}(\theta_t)\right\}\dd t.
\end{equation}

We pass to $g$ in each term. For the endpoints, Young's inequality,
Lemma~\ref{lem:pairings} and \eqref{eq:energy} imply
\begin{equation*}
 \sup_{0\le t\le T}|F_{g_\varepsilon}(\theta_t)-F_g(\theta_t)| \le\norm{g_\varepsilon-g}_1\sup_t\norm{\theta_t}_2^2 \le C_\gamma\varepsilon^{\beta-2\gamma}\norm{\theta_0}_2^2.
\end{equation*}
For $\mathcal B_ag_\varepsilon$ and $\mathcal B_vg_\varepsilon$,
the Schur bound in
Definition~\ref{def:kernel} gives
\begin{equation}\label{eq:energy-B-error}
 \begin{aligned}
 &\int_0^S\left\{\sum_{a=0}^J
 |F_{\mathcal B_ag_\varepsilon}(\theta_t)-F_{\mathcal B_ag}(\theta_t)|
                  +|F_{\mathcal B_vg_\varepsilon}(\theta_t)-F_{\mathcal B_vg}(\theta_t)|\right\}\dd t\\
 &\le C_\gamma\varepsilon^{\beta-2\gamma}
       \left(\sum_{a=0}^J\lambda_a+\norm{Dv}_\infty\right)
                                  \int_0^S\norm{\theta_t}_2^2\dd t \le C_\gamma\varepsilon^{\beta-2\gamma}
       \left(\sum_{a=0}^N\lambda_a\right)\norm{\theta_0}_2^2.
 \end{aligned}
\end{equation}
For \eqref{eq:energy-B-error} we used
$\int_0^S\norm{\theta_t}_2^2\dd t\le T\norm{\theta_0}_2^2$,
\eqref{eq:energy} and
\begin{equation*}
\norm{Dv}_\infty\le\sum_{a=J+1}^N\norm{Du_a}_\infty
                       \le C\sum_{a=J+1}^N\lambda_a,
\end{equation*}
which follows from Lemma~\ref{lem:spatial}.
The bound in \eqref{eq:energy-B-error} tends to zero at fixed $N$; no uniformity in
$N$ is required here.

For almost every $t$, $\theta_t\in H^1$. Formula
\eqref{eq:lap-action} and Young's inequality therefore yield
\begin{equation*}
 \begin{aligned}
 &\mu\int_0^S|F_{\Delta_Yg_\varepsilon}(\theta_t)
                            -F_{\Delta_Yg}(\theta_t)|\dd t =2\mu\int_0^S\left|\langle\nabla\theta_t,
                      (g_\varepsilon-g)*\nabla\theta_t\rangle\right|\dd t\\
 &\le2\norm{g_\varepsilon-g}_1\mu\int_0^S\norm{\nabla\theta_t}_2^2\dd t
 \le\norm{g_\varepsilon-g}_1\norm{\theta_0}_2^2
 \le C_\gamma\varepsilon^{\beta-2\gamma}\norm{\theta_0}_2^2.
 \end{aligned}
\end{equation*}
The limiting quadratic expressions are integrable: the same Schur estimates
bound the integrals of $|F_{\mathcal B_ag}(\theta_t)|$ and
$|F_{\mathcal B_vg}(\theta_t)|$ by
$C_\gamma(\sum_{a=0}^N\lambda_a)\norm{\theta_0}_2^2$, while
\begin{equation*}
\mu\int_0^S|F_{\Delta_Yg}(\theta_t)|\dd t
 \le2\norm g_1\mu\int_0^S\norm{\nabla\theta_t}_2^2\dd t
 \le\norm g_1\norm{\theta_0}_2^2.
\end{equation*}
All errors are uniform in $S\in[0,T]$. Letting
$\varepsilon\downarrow0$ in \eqref{eq:regularized-energy} proves
\eqref{eq:singular-energy} for every such $S$.
\end{proof}

\section{The exact renewal identity}\label{sec:renewal}

In this section we prove the main identity leading to the proof of
Theorem~\ref{thm:main}. We apply the results of the previous section
on the renewal intervals $I_{jk}$ and then combine the resulting
identities over $[0,T]$. Recall from \eqref{eq:model-horizon} that
$T=\kappa\lceil\kappa^{-1}\rceil$, so every generation partitions
$[0,T]$ into complete renewal intervals. Since the integrand in
Theorem~\ref{thm:main} is nonnegative and $T\ge1$, an upper bound
on its integral over $[0,T]$ gives the required bound over $[0,1]$.

For each $0<\mu\le1$, choose the generation cutoff
$J=J(\mu)$ by
\begin{equation*}
 \ell_J^{1+\alpha}\le\mu<2^{1+\alpha}\ell_J^{1+\alpha},
 \qquad J\ge0,\quad0<\mu\le1.
\end{equation*}
In particular, $J=0$ when $\mu=1$. Take $N\ge J$ and a smooth
deterministic datum, and write $\theta$ for the scalar driven
by $u^{(N)}$, with $v=u^{(N)}-u^{(J)}$ as before.

We will apply Lemma~\ref{lem:product} up to three times and
integrate each identity. This produces iterated time integrals
which we shall compute explicitly, leading to certain time weights.
In order to shorten our formulas, it is convenient to introduce
these weights before stating the renewal identity.

\begin{definition}[Renewal intervals and time weights]\label{def:time-weights}
For $0\le k<T/\tau_b$, the scale-$b$ generation born at
$s_{b,k}=k\tau_b$ is active on
$I_{bk}=[s_{b,k},s_{b,k}+\tau_b)$. For $0\le t<T$, set
\begin{equation*}
 s_b(t)=\tau_b\lfloor t/\tau_b\rfloor,\qquad e_b(t)=s_b(t)+\tau_b.
\end{equation*}
Thus $s_b(t)$ and $e_b(t)$ are the birth and next renewal time of
the generation containing $t$. In particular, $e_b(t)\le T$.
For $0\le s<T$, put
\begin{equation}\label{eq:weights-h}
 h_b(s)=e_b(s)-s,
\end{equation}
\begin{equation}\label{eq:weights}
 q_{ab}(s)=\int_s^{\min(e_a(s),e_b(s))}(e_b(s)-r)\dd r.
\end{equation}
The weight $h_b(s)$ in \eqref{eq:weights-h} is the time remaining
after $s$ in the $b$-generation. The weight $q_{ab}(s)$ in
\eqref{eq:weights} is the area of the pairs of times
\begin{equation*}
s\le r\le\min(e_a(s),e_b(s)),\qquad r\le t\le e_b(s).
\end{equation*}
\end{definition}
The functions $h_b,q_{ab}$ in \eqref{eq:weights-h}--\eqref{eq:weights} satisfy
\begin{equation}\label{eq:weight-bounds}
 0\le h_b\le\frac\kappa{\lambda_b},\qquad
 0\le q_{ab}\le\tau_b\min(\tau_a,\tau_b)
       \le\frac{\kappa^2}{\lambda_b\max(\lambda_a,\lambda_b)}.
\end{equation}
The scales too coarse for a given separation are removed from the
iterated kernels, while the endpoint energy remains $F_g$.
Introduce $L\ge4$, to be fixed as $L=\kappa^{-2/\beta}$ in
Section~\ref{sec:completion}. Choose a smooth nonincreasing $\chi$
equal to one on $[0,1]$ and zero on $[2,\infty)$.

\begin{definition}[Localized kernels]
Define $g_b$ by
\begin{equation*}
 g_b'(r)=\chi(\ell_b/(Lr))g'(r),\qquad
 g_b(r)=-\int_r^{2r_0}g_b'(s)\dd s\quad(0\le r<2r_0),
\end{equation*}
and set it to zero outside $B_{2r_0}(0)$, periodically.  Define moreover
\begin{equation}\label{eq:E0-Gamma}
 E_0=\sum_{b=0}^J\mathcal B_b(g-g_b),\qquad
 \Gamma_b(z)=\frac{\sigma^2A_b^2}{2}D(z/\ell_b):D^2g_b(z),
\end{equation}
and, using the weights in \eqref{eq:weights-h}--\eqref{eq:weights},
\begin{equation}\label{eq:response-kernels}
 H_1(s)=\sum_{b=0}^Jh_b(s)\mathcal B_bg_b,\qquad
 H_2(s)=\sum_{a,b=0}^Jq_{ab}(s)\mathcal B_a\mathcal B_bg_b.
\end{equation}
\end{definition}

The localized kernels are $C^4$ and constant near zero. Differentiating
the cutoff shows, with constants independent of $L$,
\begin{equation}\label{eq:gb-derivatives}
 |D^kg_b(z)|\le C_\gamma r^{-d+\beta-2\gamma-k}
                  \ind_{\{\ell_b\le2Lr\}},\qquad1\le k\le4.
\end{equation}
For $1\le q\le4$, Leibniz's formula gives
\begin{equation*}
 g_b^{(q)}(r)=\sum_{m=0}^{q-1}\binom{q-1}{m}
 \frac{\dd^m}{\dd r^m}\left[\chi\left(\frac{\ell_b}{Lr}\right)\right]
 g^{(q-m)}(r).
\end{equation*}
The $m$th derivative of the bracket is at most $C_mr^{-m}$;
for $m>0$ it is supported where $\ell_b/(Lr)\in[1,2]$, and for
$m=0$ the bracket vanishes when $\ell_b>2Lr$. Therefore
\begin{equation*}
 |g_b^{(q)}(r)|
 \le C_\gamma\sum_{m=0}^{q-1}r^{-m}r^{-d+\beta-2\gamma-(q-m)}
                         \ind_{\{\ell_b\le2Lr\}}
 \le C_\gamma r^{-d+\beta-2\gamma-q}\ind_{\{\ell_b\le2Lr\}}.
\end{equation*}
The radial derivative formula in
\eqref{eq:g-Cartesian-derivative-estimate} gives
\begin{equation*}
 |D^qg_b(z)|\le C_q\sum_{m=1}^qr^{m-q}|g_b^{(m)}(r)|
 \le C_\gamma r^{-d+\beta-2\gamma-q}\ind_{\{\ell_b\le2Lr\}},
 \qquad1\le q\le4.
\end{equation*} For later use, $\Gamma_b$ is a regular difference kernel of zero
integral. In fact, extending it by zero outside the chart,
\begin{equation}\label{eq:Gamma-zero}
 \int_{\T^d}\Gamma_b(z)\dd z
 =\frac{\sigma^2A_b^2}{2}\int D_{im}(z/\ell_b)\partial^2_{im}g_b(z)\dd z
 =-\frac{\sigma^2A_b^2}{2}\int
        \partial_i[D_{im}(z/\ell_b)]\partial_mg_b(z)\dd z=0,
\end{equation}
In \eqref{eq:Gamma-zero} repeated indices are summed. 

\begin{proposition}[Renewal identity]\label{prop:renewal}
For the finite fields, datum, and cutoffs above, with $T$ as in
\eqref{eq:model-horizon},
\begin{equation}\label{eq:renewal}
 \begin{aligned}
 &\E F_g(\theta_T)-F_g(\theta_0)
 -\underbrace{\sum_{b,k}\tau_b^2\E F_{\Gamma_b}(\theta_{s_{b,k}})}_{=:\,\mathrm{(C)}} =\E\int_0^T\left\{F_{E_0}(\theta_t)
              +\mu F_{\Delta_Yg}(\theta_t)+F_{\mathcal B_vg}(\theta_t)\right\}\dd t\\
 &+\E\int_0^T\sum_{a,b,c=0}^Jq_{ab}(t)
                   F_{\mathcal B_c\mathcal B_a\mathcal B_bg_b}(\theta_t)\dd t +\mu\E\int_0^T\sum_{i=1}^2F_{\Delta_YH_i}(\theta_t)\dd t
       +\E\int_0^T\sum_{i=1}^2F_{\mathcal B_vH_i}(\theta_t)\dd t.
 \end{aligned}
\end{equation}
Here $0\le b\le J$ and $0\le k<T/\tau_b$.
\end{proposition}
\begin{remark}
We have highlighted the term (C) because it is the one which will lead to our estimate. In fact this is the term which will control (up to a minus sign and some additional errors) the desired Sobolev seminorm, cf. Corollary \ref{c:controllo}.
\end{remark}
\begin{proof}
Choose
$0<2\varepsilon<\min(r_0,\ell_J/(2L))$. Applying
\eqref{eq:regularized-energy} and taking expectations gives
\begin{equation}\label{eq:renewal-start}
 \begin{aligned}
 &\E F_{g_\varepsilon}(\theta_T)-F_{g_\varepsilon}(\theta_0)
 -\underbrace{\sum_{b=0}^J\int_0^T\E F_{\mathcal B_bg_b}(\theta_t)\dd t}_{=:\,\mathrm{(A)}} =\underbrace{\sum_{b=0}^J\int_0^T\E\left[
 F_{\mathcal B_bg_\varepsilon}(\theta_t)
                         -F_{\mathcal B_bg_b}(\theta_t)\right]\dd t}_{=:\,\mathrm{(R_\varepsilon)}}\\
 &+\mu\E\int_0^TF_{\Delta_Yg_\varepsilon}(\theta_t)\dd t
             +\E\int_0^TF_{\mathcal B_vg_\varepsilon}(\theta_t)\dd t.
 \end{aligned}
\end{equation}
To compute $\mathrm{(A)}$, we apply Lemma~\ref{lem:product} to
$K=\mathcal B_bg_b$, starting at $s_b(t)$, and then to
$K=\mathcal B_a\mathcal B_bg_b$, starting at the birth of scale
$\max(a,b)$. These kernels are respectively $C^3$ and
$C^2$, so both applications are permitted. We first compute their
expected values at these births, where the scalar trace is
independent of the new generation variables.

At a birth $s$ of scale $m$, let $\mathcal F_{m,s}$ be generated by
all generation data born strictly before $s$ and all families at
scales strictly less than $m$. The scalar trace $\theta_s$ is
$\mathcal F_{m,s}$-measurable, since it depends only on the drift
before $s$. The new scale-$m$ shift and coefficients are independent
of this sigma-algebra. Here $\E_{\rm new}$ means integration over
that generation's common shift and all its coefficients, with the
older and coarser data fixed. For the mixed terms below we first
average its coefficients with its shift fixed.

For a spatial function $G$ write $\delta G=G(x)-G(\widetilde x)$.
The spatial chain rules \cite[(4.16)]{FM}, applied to $g_b$, give
\begin{equation*}
 \mathcal B_bg_b=\delta u_b\cdot Dg_b,
\end{equation*}
\begin{equation*}
 \mathcal B_a\mathcal B_bg_b
 =D^2g_b[\delta u_a,\delta u_b]+Dg_b\cdot\delta(\mathcal B_au_b).
\end{equation*}
The formula for $\mathcal B_au_b$ is \cite[(4.11)]{FM}.
At the birth of $m=\max(a,b)$, if $a\ne b$, every term of the
kernel $\mathcal B_a\mathcal B_bg_b$ is linear in a new
scale-$m$ coefficient. Its
coefficient average is therefore zero. This remains true for equal
lifetimes, since the coarser families, including their same-time
births, are conditioned on. Likewise $\mathcal B_bg_b$ has zero
coefficient average at the birth of $b$.

For $a=b$, Lemma~\ref{lem:average} applies at birth, when the
centers form the shifted grid. Thus
\begin{equation*}
 \E[\mathcal B_bg_b\mid\mathcal F_{b,s}]=0,
\end{equation*}
\begin{equation*}
 \E[\mathcal B_a\mathcal B_bg_b\mid\mathcal F_{m,s}]=0\qquad(a\ne b),
\end{equation*}
\begin{equation}\label{eq:birth}
 \E[\mathcal B_b^2g_b\mid\mathcal F_{b,s}] =D^2g_b:\E_{\rm new}(\delta u_b)^{\otimes2} +Dg_b\cdot\E_{\rm new}\delta[(u_b\cdot\nabla)u_b] =2\Gamma_b.
\end{equation}
The scalar trace is used only at birth in this conditional average.
Its measurability allows us to pass the conditional expectation through the form:
\begin{equation*}
\E[F_H(\theta_s)\mid\mathcal F_{m,s}]
 =\iint\theta_s(x)\theta_s(\widetilde x)
           \E[H(x,\widetilde x)\mid\mathcal F_{m,s}]\dd x\dd\widetilde x.
\end{equation*}

Applying Lemma~\ref{lem:product} to $\mathcal B_bg_b$ and integrating
from $s_b(t)$ to $t$ now gives
\begin{equation}\label{eq:first-integration}
 \E F_{\mathcal B_bg_b}(\theta_t)
 =\sum_{a=0}^J\int_{s_b(t)}^t
             \E F_{\mathcal B_a\mathcal B_bg_b}(\theta_r)\dd r +\mu\int_{s_b(t)}^t
       \E F_{\Delta_Y\mathcal B_bg_b}(\theta_r)\dd r
       +\int_{s_b(t)}^t
       \E F_{\mathcal B_v\mathcal B_bg_b}(\theta_r)\dd r.
\end{equation}
To compute the integrals of $F_{\mathcal B_a\mathcal B_bg_b}$ in
\eqref{eq:first-integration}, apply Lemma~\ref{lem:product} to
$\mathcal B_a\mathcal B_bg_b$. With $m=\max(a,b)$, integrate
from $s_m(r)$ to $r$ and move the conditional birth value in
\eqref{eq:birth} to the left:
\begin{equation}\label{eq:second-integration}
 \begin{aligned}
 &\E F_{\mathcal B_a\mathcal B_bg_b}(\theta_r)
 -2\ind_{\{a=b\}}\E F_{\Gamma_b}(\theta_{s_b(r)})
 =\sum_{c=0}^J\int_{s_m(r)}^r
            \E F_{\mathcal B_c\mathcal B_a\mathcal B_bg_b}(\theta_s)\dd s\\
 &\quad+\mu\int_{s_m(r)}^r
            \E F_{\Delta_Y\mathcal B_a\mathcal B_bg_b}(\theta_s)\dd s
       +\int_{s_m(r)}^r
            \E F_{\mathcal B_v\mathcal B_a\mathcal B_bg_b}(\theta_s)\dd s.
 \end{aligned}
\end{equation}
The data occurring in $\mathcal B_bg_b$ are those of scale $b$.
None is replaced on $(s_b(t),t)$: a renewal at a scale $a\le b$
would also be a scale-$b$ renewal, by the nesting in
\cite[Lemma 2.1]{FM}. Similarly, the coefficients and centers in
$\mathcal B_a\mathcal B_bg_b$ have scales $a,b$; with
$m=\max(a,b)$, none is replaced on $(s_m(r),r)$.
Renewals at scales greater than $b$ or $m$, respectively, can change
$u^{(N)}$ but do not change these kernel data. At each such time the
scalar has matching $L^2$ traces, so the identities from
Lemma~\ref{lem:product} concatenate. This also covers equal
lifetimes. We leave the integrals of
$F_{\mathcal B_c\mathcal B_a\mathcal B_bg_b}$ in
\eqref{eq:second-integration} unchanged and do not apply
Lemma~\ref{lem:product} to that kernel.

Substituting \eqref{eq:second-integration} into
\eqref{eq:first-integration} and integrating in $t$ computes
$\mathrm{(A)}$ of \eqref{eq:renewal-start}. The weights in
Definition~\ref{def:time-weights}, \eqref{eq:weights-h}--\eqref{eq:weights}, arise as follows.
On a
$b$-interval $[s,e]$, where $s=s_{b,k}$ and $e=s+\tau_b$,
the diagonal term is exactly
\begin{equation*}
 2\int_s^e\int_s^t\dd r\dd t\,
         \E F_{\Gamma_b}(\theta_s)
 =2\int_s^e(t-s)\dd t\,\E F_{\Gamma_b}(\theta_s)
 =\tau_b^2\E F_{\Gamma_b}(\theta_s).
\end{equation*}
At fixed $N,J,L$, all kernels in these iterated time integrals are
bounded. Schur's estimate and \eqref{eq:energy} give
\begin{equation*}
|F_H(\theta_t)|\le\norm H_{\rm Schur}\norm{\theta_t}_2^2
                  \le C_{N,J,L,\gamma}\norm{\theta_0}_2^2
\end{equation*}
for each of them. Thus the time and probability integrals are
absolutely integrable, and Fubini's theorem applies.
For the integrals of $F_{\Delta_Y\mathcal B_bg_b}$ and
$F_{\mathcal B_v\mathcal B_bg_b}$ in \eqref{eq:first-integration},
reversing the two integrations gives
\begin{equation*}
\int_0^T\int_{s_b(t)}^t A(s)\dd s\dd t
 =\int_0^T\int_s^{e_b(s)}A(s)\dd t\dd s
 =\int_0^Th_b(s)A(s)\dd s.
\end{equation*}
For each of the three terms integrated in
\eqref{eq:second-integration}, the domain is
\begin{equation*}
0\le t\le T,\qquad s_b(t)\le r\le t,\qquad s_m(r)\le s\le r.
\end{equation*}
Up to endpoints, nesting makes this equivalent to
\begin{equation*}
0\le s\le T,\qquad
 s\le r\le\min(e_a(s),e_b(s)),\qquad
 r\le t\le e_b(s).
\end{equation*}
Indeed $s,r$ lie in the same $a$- and $b$-generations, while $r,t$
lie in the same $b$-generation. Integrating first in $t$ and $r$
therefore gives
\begin{equation*}
\int_s^{\min(e_a(s),e_b(s))}
       \int_r^{e_b(s)}\dd t\dd r=q_{ab}(s).
\end{equation*}
Using the definitions of $H_1,H_2$ in
\eqref{eq:response-kernels}, the identities just computed give
\begin{equation}\label{eq:A-minus-C}
 \begin{aligned}
 &\mathrm{(A)}-\mathrm{(C)}
 =\E\int_0^T\sum_{a,b,c=0}^Jq_{ab}(t)
 F_{\mathcal B_c\mathcal B_a\mathcal B_bg_b}(\theta_t)\dd t\\
 &\quad+\mu\E\int_0^T\sum_{i=1}^2F_{\Delta_YH_i}(\theta_t)\dd t
       +\E\int_0^T\sum_{i=1}^2F_{\mathcal B_vH_i}(\theta_t)\dd t.
 \end{aligned}
\end{equation}
Thus substitution in \eqref{eq:renewal-start} puts
$\E F_{g_\varepsilon}(\theta_T)-F_{g_\varepsilon}(\theta_0)
-\mathrm{(C)}$ on the left, as in \eqref{eq:renewal}.

Finally let $\varepsilon\downarrow0$ in \eqref{eq:renewal-start}.
Our choice $2\varepsilon<\ell_J/(2L)$ puts the changed region
inside the part where every $g_b$ is constant. Hence all its
iterated kernels and birth averages are unchanged. Corollary~\ref{cor:singular-product}
gives convergence of the endpoint values and the integrals of
$F_{\Delta_Yg_\varepsilon}$ and $F_{\mathcal B_vg_\varepsilon}$
at fixed $N$. For $\mathrm{(R_\varepsilon)}$ in
\eqref{eq:renewal-start},
Lemma~\ref{lem:pairings} gives
\begin{equation*}
\sum_{b=0}^J\left[
 F_{\mathcal B_bg_\varepsilon}(\theta_t)
                         -F_{\mathcal B_bg_b}(\theta_t)\right]
 \longrightarrow F_{E_0}(\theta_t)
\end{equation*}
in integrated expectation. Hence $\mathrm{(R_\varepsilon)}$
converges to $\E\int_0^T F_{E_0}(\theta_t)\dd t$. Combined
with \eqref{eq:A-minus-C}, this proves
\eqref{eq:renewal}. No time derivative of the velocity is used.
\end{proof}

\section{The principal contribution}\label{sec:principal}
Recall the term in \eqref{eq:renewal}:
\begin{equation}\label{eq:covariance-recall}
 \mathrm{(C)}=\sum_{b=0}^J\sum_{0\le k<T/\tau_b}
                 \tau_b^2\E F_{\Gamma_b}(\theta_{s_{b,k}}).
\end{equation}
The goal of this section is to prove \eqref{eq:covariance-control}
in Corollary~\ref{c:controllo}. We first compare
\eqref{eq:covariance-recall} with the integral of
$\sum_b\tau_bF_{\Gamma_b}(\theta_t)$. The covariance in $\mathrm{(C)}$
is evaluated at the births of the generations, whereas coercivity
requires the scalar at the current time. After this comparison,
we combine the current covariance with $\mu F_{\Delta_Yg}(\theta_t)$;
the latter compensates for the covariance of the omitted scales $b>J$.

\begin{lemma}[Sampling]\label{lem:sampling}
For $T$ as in \eqref{eq:model-horizon}, uniformly in $N\ge J$
and $0<\mu\le1$,
\begin{equation}\label{eq:sampling}
 \begin{aligned}
 &-\E\int_0^T\sum_{b=0}^J\tau_bF_{\Gamma_b}(\theta_t)\dd t
 \le-\mathrm{(C)}
 +C_\gamma\kappa^2(1+\log L)
                     \E\int_0^T\norm{\theta_t}_{\dot H^\gamma}^2\dd t\\
 &\quad+C_\gamma\kappa^2(1+\log L)\mu
                         \E\int_0^T\norm{\nabla\theta_t}_2^2\dd t,
 \end{aligned}
\end{equation}
where $\mathrm{(C)}$ is defined in \eqref{eq:renewal}.
\end{lemma}

\begin{lemma}[Covariance and diffusion]\label{lem:coercivity}
For sufficiently small $\kappa$ depending on $\gamma$ and the model,
and every real $f\in H^1(\T^d)$,
\begin{equation}\label{eq:coercive}
 c_\gamma\kappa\norm f_{\dot H^\gamma}^2 \le-\sum_{b=0}^J\tau_bF_{\Gamma_b}(f)-\mu F_{\Delta_Yg}(f) +C_\gamma\kappa L^{-\beta}\norm f_{\dot H^\gamma}^2 +C_\gamma\norm f_2^2.
\end{equation}
The constants are independent of $L,J,N,\mu$.
\end{lemma}

\begin{corollary}\label{c:controllo}
Under the hypotheses of Lemmas~\ref{lem:sampling} and
\ref{lem:coercivity},
\begin{equation}\label{eq:covariance-control}
 \begin{aligned}
 &\E\int_0^T\norm{\theta_t}_{\dot H^\gamma}^2\dd t
 \le-\frac{\mathrm{(C)}}{c_\gamma\kappa}
       -\frac{\mu}{c_\gamma\kappa}\E\int_0^TF_{\Delta_Yg}(\theta_t)\dd t +C_\gamma[L^{-\beta}+\kappa(1+\log L)]
                       \E\int_0^T\norm{\theta_t}_{\dot H^\gamma}^2\dd t\\
 &+\frac{C_\gamma}{\kappa}\E\int_0^T\norm{\theta_t}_2^2\dd t
       +C_\gamma\kappa(1+\log L)\mu
                         \E\int_0^T\norm{\nabla\theta_t}_2^2\dd t.
 \end{aligned}
\end{equation}
Here $c_\gamma$ is the constant in \eqref{eq:coercive}; the
constants are independent of $L,J,N,\mu,\kappa$.
\end{corollary}

\begin{proof}[Proof of Lemma~\ref{lem:sampling}]
The term $\mathrm{(C)}$ is the sum in \eqref{eq:covariance-recall},
where $s_{b,k}=k\tau_b$. The quantity with which we compare it is
\begin{equation*}
 \E\int_0^T\sum_{b=0}^J\tau_bF_{\Gamma_b}(\theta_t)\dd t
 =\sum_{b=0}^J\sum_{0\le k<T/\tau_b}\tau_b
       \int_{s_{b,k}}^{s_{b,k}+\tau_b}\E F_{\Gamma_b}(\theta_t)\dd t.
\end{equation*}
Since every interval has length $\tau_b$, their difference satisfies
\begin{equation}\label{eq:sampling-difference}
 \left|\mathrm{(C)}-\E\int_0^T\sum_{b=0}^J \tau_bF_{\Gamma_b}(\theta_t)\dd t\right| =\left|\sum_{b=0}^J\sum_{0\le k<T/\tau_b}\tau_b \int_{s_{b,k}}^{s_{b,k}+\tau_b} \E\left[F_{\Gamma_b}(\theta_{s_{b,k}}) -F_{\Gamma_b}(\theta_t)\right]\dd t\right|.
\end{equation}
We first estimate the derivatives of $F_{\Gamma_b}(\theta_t)$.

For $k=1,2$ we claim
\begin{equation}\label{eq:Gamma-derivatives}
 \sum_{b=0}^J\tau_b^2|D^k\Gamma_b(z)|
       \le C_\gamma\kappa^2(1+\log L)
                      r^{-d+\beta-2\gamma-k}.
\end{equation}
For $0\le h\le2$, evenness, \eqref{eq:D-bounds}, and constancy
outside $B_{1/4}(0)$ give
\begin{equation*}
|D_z^h[D(z/\ell_b)]|
       \le C\min(r^2/\ell_b^2,1)r^{-h}.
\end{equation*}
Differentiation and \eqref{eq:gb-derivatives} therefore bound the
left side of \eqref{eq:Gamma-derivatives} by
\begin{equation}\label{eq:Gamma-derivative-sum}
C_\gamma\kappa^2 r^{-d+\beta-2\gamma-2-k}
 \sum_{\ell_b\le2Lr}\ell_b^2\min(r^2/\ell_b^2,1).
\end{equation}
The sum in \eqref{eq:Gamma-derivative-sum} satisfies
\begin{equation*}
 \begin{aligned}
 &\sum_{\ell_b\le2Lr}\ell_b^2\min(r^2/\ell_b^2,1)
 \le\sum_{\ell_b<r}\ell_b^2
       +r^2\#\{b:r\le\ell_b\le2Lr\}\\
 &\le Cr^2+r^2(1+\log_2(2L))
 \le C(1+\log L)r^2.
 \end{aligned}
\end{equation*}
Substitution in \eqref{eq:Gamma-derivative-sum} proves
\eqref{eq:Gamma-derivatives} with an $L$-independent constant.

Applying Lemma~\ref{lem:product} to the deterministic difference
kernel $\Gamma_b$ gives
\begin{equation}\label{eq:Gamma-evolution}
\frac\dd{\dd t}F_{\Gamma_b}(\theta_t)
 =F_{[u^{(N)}(x)-u^{(N)}(\widetilde x)]\cdot D\Gamma_b}(\theta_t)
                         +2\mu F_{\Delta_z\Gamma_b}(\theta_t).
\end{equation}
The function $t\mapsto F_{\Gamma_b}(\theta_t)$ is absolutely
continuous on $[0,T]$: $\Gamma_b$ is deterministic and regular,
and the continuous $L^2$ trace of $\theta$ makes its values match
at every renewal, as in Corollary~\ref{cor:singular-product}.
By \eqref{eq:Gamma-zero} and
Lemma~\ref{lem:pairings}(\ref{item:cancellation}), both kernels
in \eqref{eq:Gamma-evolution} satisfy \eqref{eq:increments}.
Using \eqref{eq:Gamma-derivatives} with $k=1$ and the spatial
bound \eqref{eq:intro-regularity}, we obtain
\begin{equation}\label{eq:Gamma-u-bound}
 \begin{aligned}
 &\sum_{b=0}^J\tau_b^2
 \left|F_{[u^{(N)}(x)-u^{(N)}(\widetilde x)]\cdot D\Gamma_b}(\theta_t)\right| \le C_\gamma\kappa^2(1+\log L)
 \int_{|z|<2r_0}\norm{\theta_t(\cdot+z)-\theta_t}_2^2
                    r^{\alpha-d+\beta-2\gamma-1}\dd z\\
 &\le C_\gamma\kappa^2(1+\log L)
                              \norm{\theta_t}_{\dot H^\gamma}^2.
 \end{aligned}
\end{equation}
Here $\alpha+\beta=1$, and we used \eqref{eq:fractional-upper}.
Likewise, \eqref{eq:Gamma-derivatives} with $k=2$ and
\eqref{eq:H1-translation} give
\begin{equation}\label{eq:Gamma-lap-bound}
 \begin{aligned}
 &2\mu\sum_{b=0}^J\tau_b^2|F_{\Delta_z\Gamma_b}(\theta_t)|
 \le C_\gamma\mu\kappa^2(1+\log L)\norm{\nabla\theta_t}_2^2
                  \int_0^{2r_0}r^{\beta-2\gamma-1}\dd r\\
 &\le C_\gamma\mu\kappa^2(1+\log L)\norm{\nabla\theta_t}_2^2.
 \end{aligned}
\end{equation}
The integral is finite since $\beta-2\gamma>0$. Thus
\eqref{eq:Gamma-evolution}, \eqref{eq:Gamma-u-bound} and
\eqref{eq:Gamma-lap-bound} imply, almost everywhere,
\begin{equation}\label{eq:Gamma-time-derivative-bound}
 \sum_{b=0}^J\tau_b^2\left|\frac\dd{\dd t}F_{\Gamma_b}(\theta_t)\right|
 \le C_\gamma\kappa^2(1+\log L)
       \left(\norm{\theta_t}_{\dot H^\gamma}^2
                         +\mu\norm{\nabla\theta_t}_2^2\right).
\end{equation}
At fixed $N,\mu$, this is integrable in time and probability by
\eqref{eq:energy} and $\gamma<1$.

On an interval $[s,e]$, where $s=s_{b,k}$ and $e=s+\tau_b$,
$F_{\Gamma_b}(\theta_s)$ is constant in $t$. Hence
\begin{equation*}
 \frac\dd{\dd t}\left[F_{\Gamma_b}(\theta_s)-F_{\Gamma_b}(\theta_t)\right]
 =-\frac\dd{\dd t}F_{\Gamma_b}(\theta_t)\qquad(s<t<e).
\end{equation*}
The fundamental theorem of calculus gives
\begin{equation}\label{eq:sampling-FTC}
 F_{\Gamma_b}(\theta_s)-F_{\Gamma_b}(\theta_t)
 =-\int_s^t\frac\dd{\dd r}F_{\Gamma_b}(\theta_r)\dd r.
\end{equation}
Integrating \eqref{eq:sampling-FTC} in $t$ and reversing the order
of integration yields
\begin{equation*}
 \tau_b^2 F_{\Gamma_b}(\theta_s) -\tau_b\int_s^eF_{\Gamma_b}(\theta_t)\dd t =-\tau_b\int_s^e\int_s^t \frac\dd{\dd r}F_{\Gamma_b}(\theta_r)\dd r\dd t =-\tau_b\int_s^e(e-r)\frac\dd{\dd r} F_{\Gamma_b}(\theta_r)\dd r.
\end{equation*}
In particular, \eqref{eq:sampling-difference} and
\eqref{eq:sampling-FTC}, followed by Tonelli's theorem, give
\begin{equation*}
 \begin{aligned}
 &\left|\mathrm{(C)}-\E\int_0^T\sum_{b=0}^J
                           \tau_bF_{\Gamma_b}(\theta_t)\dd t\right| \le\sum_{b,k}\tau_b\int_{s_{b,k}}^{s_{b,k}+\tau_b}
            \int_{s_{b,k}}^t\E\left|\frac\dd{\dd r}
                                F_{\Gamma_b}(\theta_r)\right|\dd r\dd t\\
 &=\sum_{b,k}\tau_b\int_{s_{b,k}}^{s_{b,k}+\tau_b}
       (s_{b,k}+\tau_b-r)\E\left|\frac\dd{\dd r}
                                F_{\Gamma_b}(\theta_r)\right|\dd r \le\E\int_0^T\sum_{b=0}^J\tau_b^2
                   \left|\frac\dd{\dd r}F_{\Gamma_b}(\theta_r)\right|\dd r.
 \end{aligned}
\end{equation*}
We used $0\le s_{b,k}+\tau_b-r\le\tau_b$ and the partition of
$[0,T]$ into the scale-$b$ intervals. Substitution of
\eqref{eq:Gamma-time-derivative-bound} proves the absolute-value
estimate
\begin{equation*}
 \begin{aligned}
 &\left|\mathrm{(C)}-\E\int_0^T\sum_{b=0}^J
                         \tau_bF_{\Gamma_b}(\theta_t)\dd t\right|\\
 &\quad\le C_\gamma\kappa^2(1+\log L)
                     \E\int_0^T\norm{\theta_t}_{\dot H^\gamma}^2\dd t
       +C_\gamma\kappa^2(1+\log L)\mu
                         \E\int_0^T\norm{\nabla\theta_t}_2^2\dd t.
 \end{aligned}
\end{equation*}
The inequality \eqref{eq:sampling} of
Lemma~\ref{lem:sampling} follows immediately.
\end{proof}

\begin{proof}[Proof of Lemma~\ref{lem:coercivity}]
Let $\mathcal D_{\le J}$ denote \eqref{eq:averaged-covariance} truncated at $J$. Then
\begin{equation*}
 \norm{\mathcal D-\mathcal D_{\le J}}_{\rm op}
 \le C\kappa\sum_{b>J}\ell_b^{2-\beta}
 \le C\kappa\ell_J^{1+\alpha}\le C\kappa\mu.
\end{equation*}
On $0<r<r_0$, Lemma~\ref{lem:g} therefore gives
\begin{equation}\label{eq:compensation}
 (\mathcal D_{\le J}+2\mu\Id):D^2g
 \ge\kappa r^{-d-2\gamma}+(2\mu-C_\gamma\kappa\mu)\Delta g.
\end{equation}
Choose $\kappa$ so the coefficient of $\Delta g$ is nonnegative.
This is where direct diffusion compensates for the omitted fine
covariance.

We estimate $D^2g_b-D^2g$. With $e=z/r$,
\begin{equation*}
D^2g_b-D^2g=[\chi(\ell_b/(Lr))-1]D^2g
       -\frac{\ell_b}{Lr^2}\chi'(\ell_b/(Lr))g' e\otimes e.
\end{equation*}
It vanishes unless $\ell_b>Lr$, and its magnitude is at most
$C_\gamma r^{-d+\beta-2\gamma-2}$. Consequently
\begin{equation}\label{eq:localization}
 \left|\frac{\sigma^2}{2}\sum_{b=0}^J\tau_bA_b^2 D(z/\ell_b):(D^2g_b-D^2g)\right| \le C_\gamma r^{-d+\beta-2\gamma-2} \kappa r^2\sum_{\ell_b>Lr}\ell_b^{-\beta} \le C_\gamma\kappa L^{-\beta}r^{-d-2\gamma}.
\end{equation}
By Lemma~\ref{lem:pairings}(\ref{item:cancellation})--(\ref{item:lap-increments}), the sum
$\mu F_{\Delta_Yg}(f)+\sum_{b=0}^J\tau_bF_{\Gamma_b}(f)$ equals
\begin{equation}\label{eq:coercivity-increments}
-\frac12\int_{\T^d}\norm{f(\cdot+z)-f}_2^2
              \left[2\mu\Delta g(z)+\sum_{b=0}^J\tau_b\Gamma_b(z)\right]\dd z.
\end{equation}
Equations~\eqref{eq:compensation} and \eqref{eq:localization}
bound the integral in \eqref{eq:coercivity-increments} over
$B_{r_0}(0)$ by
\begin{equation}\label{eq:coercivity-inner-bound}
-\frac\kappa2\int_{|z|<r_0}
          \frac{\norm{f(\cdot+z)-f}_2^2}{r^{d+2\gamma}}\dd z
 +C_\gamma\kappa L^{-\beta}\int_{|z|<r_0}
          \frac{\norm{f(\cdot+z)-f}_2^2}{r^{d+2\gamma}}\dd z.
\end{equation}
On $B_{2r_0}(0)\setminus B_{r_0}(0)$, all second derivatives of $g,g_b$ are bounded
uniformly in $L,b$, while $\sum_b\tau_bA_b^2\le C\kappa$ and
$\mu\le1$. Thus the integral in \eqref{eq:coercivity-increments} over
$r_0\le|z|\le2r_0$ is bounded in absolute value by
$C_\gamma\norm f_2^2$.
Using \eqref{eq:fractional} for the negative integral in
\eqref{eq:coercivity-inner-bound}, and \eqref{eq:fractional-upper}
for the positive integral, gives
\begin{equation*}
 \mu F_{\Delta_Yg}(f)+\sum_b\tau_bF_{\Gamma_b}(f)
 \le-c_\gamma\kappa\norm f_{\dot H^\gamma}^2
       +C_\gamma\kappa L^{-\beta}\norm f_{\dot H^\gamma}^2
       +C_\gamma\norm f_2^2.
\end{equation*}
Moving the first term on the right to the left yields
\begin{equation*}
 c_\gamma\kappa\norm f_{\dot H^\gamma}^2
 \le-\mu F_{\Delta_Yg}(f)-\sum_b\tau_bF_{\Gamma_b}(f)
       +C_\gamma\kappa L^{-\beta}\norm f_{\dot H^\gamma}^2
       +C_\gamma\norm f_2^2.
\end{equation*}
\end{proof}

\begin{proof}[Proof of Corollary~\ref{c:controllo}]
The scalar energy identity \eqref{eq:energy} gives
\begin{equation*}
 \E\int_0^T\norm{\theta_t}_{H^1}^2\dd t
 \le\left(T+\frac1{2\mu}\right)\norm{\theta_0}_2^2<\infty.
\end{equation*}
In particular, $\theta_t\in H^1$ almost everywhere. Since $\gamma<1$,
\begin{equation*}
 \E\int_0^T\norm{\theta_t}_{\dot H^\gamma}^2\dd t
 \le\E\int_0^T\norm{\theta_t}_{H^1}^2\dd t<\infty.
\end{equation*}
Lemma~\ref{lem:pairings}(\ref{item:lap-convergence}) and
\eqref{eq:energy} also give
\begin{equation*}
 \mu\E\int_0^T|F_{\Delta_Yg}(\theta_t)|\dd t \le2\norm g_1\mu\E\int_0^T\norm{\nabla\theta_t}_2^2\dd t \le\norm g_1\norm{\theta_0}_2^2.
\end{equation*}
At fixed $J,L$, each $\Gamma_b$ is regular and integrable, so
\begin{equation*}
 \E\int_0^T\sum_{b=0}^J\tau_b|F_{\Gamma_b}(\theta_t)|\dd t
 \le T\sum_{b=0}^J\tau_b\norm{\Gamma_b}_1\norm{\theta_0}_2^2<\infty.
\end{equation*}
We can therefore set $f=\theta_t$ in \eqref{eq:coercive}, integrate
in $t$ and take expectations:
\begin{equation*}
 \begin{aligned}
 &c_\gamma\kappa\E\int_0^T\norm{\theta_t}_{\dot H^\gamma}^2\dd t
 \le-\E\int_0^T\sum_{b=0}^J\tau_bF_{\Gamma_b}(\theta_t)\dd t
       -\mu\E\int_0^TF_{\Delta_Yg}(\theta_t)\dd t\\
 &\quad+C_\gamma\kappa L^{-\beta}
               \E\int_0^T\norm{\theta_t}_{\dot H^\gamma}^2\dd t
       +C_\gamma\E\int_0^T\norm{\theta_t}_2^2\dd t.
 \end{aligned}
\end{equation*}
Substituting \eqref{eq:sampling} for
$-\E\int_0^T\sum_b\tau_bF_{\Gamma_b}(\theta_t)\dd t$ yields
\begin{equation*}
 \begin{aligned}
 &c_\gamma\kappa\E\int_0^T\norm{\theta_t}_{\dot H^\gamma}^2\dd t
 \le-\mathrm{(C)}-\mu\E\int_0^TF_{\Delta_Yg}(\theta_t)\dd t +C_\gamma[\kappa L^{-\beta}+\kappa^2(1+\log L)]
                 \E\int_0^T\norm{\theta_t}_{\dot H^\gamma}^2\dd t\\
 &+C_\gamma\E\int_0^T\norm{\theta_t}_2^2\dd t
       +C_\gamma\kappa^2(1+\log L)\mu
                 \E\int_0^T\norm{\nabla\theta_t}_2^2\dd t.
 \end{aligned}
\end{equation*}
Dividing by $c_\gamma\kappa$ and enlarging $C_\gamma$ gives
\begin{equation*}
 \begin{aligned}
 &\E\int_0^T\norm{\theta_t}_{\dot H^\gamma}^2\dd t
 \le-\frac{\mathrm{(C)}}{c_\gamma\kappa}
       -\frac{\mu}{c_\gamma\kappa}\E\int_0^TF_{\Delta_Yg}(\theta_t)\dd t +C_\gamma
       [L^{-\beta}+\kappa(1+\log L)]
                 \E\int_0^T\norm{\theta_t}_{\dot H^\gamma}^2\dd t\\
 &+\frac{C_\gamma}{\kappa}\E\int_0^T\norm{\theta_t}_2^2\dd t
       +C_\gamma\kappa(1+\log L)\mu
                 \E\int_0^T\norm{\nabla\theta_t}_2^2\dd t.
 \end{aligned}
\end{equation*}
\end{proof}

\section{Estimates for the remaining operators}\label{sec:estimates}
We bound the integrals of $F_{E_0}$ and
$\sum_{a,b,c}q_{ab}F_{\mathcal B_c\mathcal B_a\mathcal B_bg_b}$
in \eqref{eq:renewal} by a small multiple of the Sobolev integral.
For $\mu F_{\Delta_YH_i}$ and $F_{\mathcal B_vH_i}$ we also use
$\mu\int_0^T\norm{\nabla\theta_t}_2^2\dd t$ from \eqref{eq:energy}.
Only the coefficients in the former bounds must be uniform in $L$.

\subsection{Estimates for the spatial derivatives}
We use the material-derivative estimates of
\cite[Lemma 4.3]{FM}, with the configuration convention of their
Definition 4.2 and Remark 4.4.
The expressions $\mathcal B_au_b(y)$ are defined by
\eqref{eq:B-ub}, with $y$ and the centers independent. Their
values at $x,\widetilde x$ are taken only after applying the operator.

\begin{lemma}
At configurations satisfying the separation in Lemma~\ref{lem:spatial},
\begin{equation}\label{eq:material}
 \norm{\mathcal B_au_b}_\infty\le C\lambda_{\min(a,b)}A_{\max(a,b)},
\end{equation}
\begin{equation}\label{eq:material-gradient}
 \norm{D_y(\mathcal B_au_b)}_\infty\le C\lambda_a\lambda_b.
\end{equation}
For $0<r=|x-\widetilde x|\le1/4$ and $i\le j\le k$ the sorted indices
$a,b,c$, one has
\begin{equation}\label{eq:velocity-increment}
 |\delta u_b|\le Cr\min(\lambda_b,A_b/r),
\end{equation}
\begin{equation}\label{eq:material-increment}
 |\delta(\mathcal B_au_b)|\le Cr\lambda_{\min(a,b)}\min(\lambda_{\max(a,b)},A_{\max(a,b)}/r),
\end{equation}
\begin{equation}\label{eq:material-increments}
 |\delta(\mathcal B_c\mathcal B_au_b)|\le Cr\lambda_i\lambda_j\min(\lambda_k,A_k/r).
\end{equation}
All constants are independent of the indices and cutoffs.
\end{lemma}
The amplitude bound is \cite[(4.13)]{FM}, the bound for $D_y(\mathcal B_au_b)$
is proved in \cite[proof of Lemma 5.4, Step 2, p. 45]{FM}, and the
increment bounds are \cite[Lemma 4.3, (4.6)--(4.8)]{FM}. Their
hypotheses are the bounded coefficients, separated centers and
$0<r\le1/4$; all hold here.

\begin{lemma}\label{lem:transport}
For real $f\in H^\gamma(\T^d)$, at each time outside renewal endpoints,
\begin{equation*}
 |F_{E_0}(f)|+
 \left|\sum_{a,b,c=0}^Jq_{ab}
                    F_{\mathcal B_c\mathcal B_a\mathcal B_bg_b}(f)\right|
 \le C_\gamma\left[L^{-\beta}+\kappa^2(1+\log L)^2\right]
                                      \norm f_{\dot H^\gamma}^2.
\end{equation*}
The constant is independent of $L$.
\end{lemma}
\begin{proof}
For this calculation only, put
$\rho_a(r)=\min(\lambda_a,A_a/r)$.
For each triple $(a,b,c)$, let $(i,j,k)$ be its nondecreasing
rearrangement, with repetitions retained.
The indices $i,j$ specify the two factors estimated by their spatial
derivative bounds $\lambda_i,\lambda_j$. At the largest index $k$
we retain both bounds through $\rho_k(r)=\min(\lambda_k,A_k/r)$.
The operators remain in the order
$\mathcal B_c\mathcal B_a\mathcal B_b$, and $g_b$ retains its
original index $b$.
The identity \cite[(4.17)]{FM}, applied to $g_b$, reads
\begin{equation}\label{eq:cubic-chain-expanded}
 \begin{aligned}
 &\mathcal B_c\mathcal B_a\mathcal B_bg_b
 =D^3g_b[\delta u_c,\delta u_a,\delta u_b] +D^2g_b[\delta(\mathcal B_cu_a),\delta u_b]
       +D^2g_b[\delta u_a,\delta(\mathcal B_cu_b)]\\
 &+D^2g_b[\delta u_c,\delta(\mathcal B_au_b)]
       +Dg_b\cdot\delta(\mathcal B_c\mathcal B_au_b).
 \end{aligned}
\end{equation}
We estimate each expression using
\eqref{eq:gb-derivatives} and
\eqref{eq:velocity-increment}--\eqref{eq:material-increments}.
Since $\rho_i\le\lambda_i$ and $\rho_j\le\lambda_j$,
\begin{equation}\label{eq:cubic-three-increments}
 |D^3g_b[\delta u_c,\delta u_a,\delta u_b]| \le C_\gamma r^{-d+\beta-2\gamma-3} r^3\rho_i\rho_j\rho_k\ind_{\{\ell_b\le2Lr\}} \le C_\gamma r^{-d+\beta-2\gamma} \lambda_i\lambda_j\rho_k\ind_{\{\ell_b\le2Lr\}}.
\end{equation}
For any permutation $(p,q,h)$ of $(a,b,c)$ we have
\begin{equation}\label{eq:cubic-pair-coefficient}
 \lambda_{\min(p,q)}\rho_{\max(p,q)}\rho_h
                  \le\lambda_i\lambda_j\rho_k.
\end{equation}
Indeed, if $h=k$, then
\begin{equation*}
 \lambda_{\min(p,q)}\rho_{\max(p,q)}\rho_h
 =\lambda_i\rho_j\rho_k\le\lambda_i\lambda_j\rho_k.
\end{equation*}
If $h<k$, then $\max(p,q)=k$ and
$\{\min(p,q),h\}=\{i,j\}$, with multiplicities. Thus
\begin{equation*}
 \lambda_{\min(p,q)}\rho_{\max(p,q)}\rho_h
 \le\lambda_{\min(p,q)}\lambda_h\rho_k
 =\lambda_i\lambda_j\rho_k.
\end{equation*}
Using \eqref{eq:material-increment} and
\eqref{eq:cubic-pair-coefficient}, we obtain
\begin{equation}\label{eq:cubic-two-increments}
 \begin{aligned}
 &|D^2g_b[\delta(\mathcal B_pu_q),\delta u_h]|
 \le C_\gamma r^{-d+\beta-2\gamma-2}
       r^2\lambda_{\min(p,q)}\rho_{\max(p,q)}\rho_h
                              \ind_{\{\ell_b\le2Lr\}}\\
 &\le C_\gamma r^{-d+\beta-2\gamma}
       \lambda_i\lambda_j\rho_k\ind_{\{\ell_b\le2Lr\}}.
 \end{aligned}
\end{equation}
The three expressions involving $D^2g_b$ in
\eqref{eq:cubic-chain-expanded} correspond respectively to
$(p,q,h)=(c,a,b)$, $(c,b,a)$, and $(a,b,c)$; symmetry of $D^2g_b$ allows
us to exchange its two vector arguments. Finally,
\eqref{eq:material-increments} gives
\begin{equation}\label{eq:cubic-one-increment}
 |Dg_b\cdot\delta(\mathcal B_c\mathcal B_au_b)| \le C_\gamma r^{-d+\beta-2\gamma-1} r\lambda_i\lambda_j\rho_k\ind_{\{\ell_b\le2Lr\}} =C_\gamma r^{-d+\beta-2\gamma} \lambda_i\lambda_j\rho_k\ind_{\{\ell_b\le2Lr\}}.
\end{equation}
Adding \eqref{eq:cubic-three-increments}, the three instances of
\eqref{eq:cubic-two-increments}, and \eqref{eq:cubic-one-increment}
in \eqref{eq:cubic-chain-expanded}, we conclude that
\begin{equation}\label{eq:third-pointwise}
 |\mathcal B_c\mathcal B_a\mathcal B_bg_b|
 \le C_\gamma r^{-d+\beta-2\gamma}\lambda_i\lambda_j\rho_k(r)
                                      \ind_{\{\ell_b\le2Lr\}}.
\end{equation}

The scale calculation in \cite[proof of Lemma 5.1, fourth term
of (4.20), pp. 42--43, using (5.10)]{FM} gives
\begin{equation}\label{eq:cubic-scale-sum}
\sum_{a,b,c}q_{ab}\lambda_i\lambda_j\rho_k(r)
           \ind_{\{\ell_b\le2Lr\}}
 \le C\kappa^2r^{-\beta}(1+\log L)^2.
\end{equation}
For these sorted indices, the
coefficient in the cited scale calculation is
\begin{equation*}
 \lambda_{\min(a,b)}\lambda_{\min(c,\max(a,b))}
 \rho_{\max(a,b,c)}=\lambda_i\lambda_j\rho_k.
\end{equation*}
In \eqref{eq:cubic-scale-sum}, $i,j,k$ are determined anew by
each summand $(a,b,c)$; the sum is still over the original indices
$a,b,c$, with the factors $q_{ab}$ and
$\ind_{\{\ell_b\le2Lr\}}$ unchanged.
The calculation uses only the upper bound in
\eqref{eq:weight-bounds}. Its localized index range contains ours,
since $2Lr\le3L\max(r,\ell_J)$, and
$\max(r,\ell_J)^{-\beta}\le r^{-\beta}$. Thus
\eqref{eq:third-pointwise} yields
\begin{equation}\label{eq:cubic-kernel-bound}
\sum_{a,b,c}q_{ab}|\mathcal B_c\mathcal B_a\mathcal B_bg_b|
 \le C_\gamma\kappa^2(1+\log L)^2r^{-d-2\gamma}.
\end{equation}
For $E_0$, only $\ell_b>Lr$ occur, and
\begin{equation}\label{eq:E0-kernel-bound}
|E_0|\le C_\gamma r^{-d+\beta-2\gamma}
                \sum_{\ell_b>Lr}\lambda_b
       \le C_\gamma L^{-\beta}r^{-d-2\gamma}.
\end{equation}
By Lemma~\ref{lem:pairings}(\ref{item:cancellation}), the kernels
$\mathcal B_c\mathcal B_a\mathcal B_bg_b$ have zero row and column
integrals. The same holds for the integrable
kernels $\mathcal B_b(g-g_b)$ defining $E_0$ in
\eqref{eq:E0-Gamma}, by their Schur limits. The increment identity
\eqref{eq:increments}, the bound \eqref{eq:E0-kernel-bound}, and
\eqref{eq:fractional-upper} give
\begin{equation}\label{eq:E0-quadratic-estimate}
 |F_{E_0}(f)|
 \le C_\gamma L^{-\beta}\int_{|z|<2r_0}
        \frac{\norm{f(\cdot+z)-f}_2^2}{|z|^{d+2\gamma}}\dd z
 \le C_\gamma L^{-\beta}\norm f_{\dot H^\gamma}^2.
\end{equation}
Likewise, \eqref{eq:cubic-kernel-bound} gives
\begin{equation}\label{eq:cubic-quadratic-estimate}
 \begin{aligned}
 &\left|\sum_{a,b,c}q_{ab}F_{\mathcal B_c\mathcal B_a\mathcal B_bg_b}(f)\right|
 \le\frac12\iint|f(x)-f(\widetilde x)|^2
       \sum_{a,b,c}q_{ab}|\mathcal B_c\mathcal B_a\mathcal B_bg_b|
                                            \dd x\dd\widetilde x\\
 &\le C_\gamma\kappa^2(1+\log L)^2
       \int_{|z|<2r_0}\frac{\norm{f(\cdot+z)-f}_2^2}{|z|^{d+2\gamma}}\dd z \le C_\gamma\kappa^2(1+\log L)^2\norm f_{\dot H^\gamma}^2.
 \end{aligned}
\end{equation}
Adding \eqref{eq:E0-quadratic-estimate} and
\eqref{eq:cubic-quadratic-estimate} gives
\begin{equation*}
 |F_{E_0}(f)|+
 \left|\sum_{a,b,c}q_{ab}F_{\mathcal B_c\mathcal B_a\mathcal B_bg_b}(f)\right|
 \le C_\gamma[L^{-\beta}+\kappa^2(1+\log L)^2]\norm f_{\dot H^\gamma}^2.
\end{equation*}
\end{proof}

\subsection{The operators \texorpdfstring{$H_1$ and $H_2$}{H1 and H2}}
We first estimate the linear operators $H_1,H_2$ defined in
\eqref{eq:response-kernels}. Their $H^1$ bounds will control the
bilinear operators $\mu F_{\Delta_YH_i}$ and $F_{\mathcal B_vH_i}$
in \eqref{eq:renewal}.

\begin{lemma}
For $f\in H^\gamma(\T^d)$ and $i=1,2$,
\begin{equation*}
 \norm{H_if}_{H^1}\le C_{\gamma,L}\kappa^i
                  \ell_J^{-\alpha-\gamma}\norm f_{\dot H^\gamma}.
\end{equation*}
Both operators annihilate constants and have mean-zero outputs.
\end{lemma}
\begin{proof}
By Lemma~\ref{lem:pairings}(\ref{item:cancellation}), every
summand in \eqref{eq:response-kernels} has zero differentiated
row integral. Hence
\begin{equation}\label{eq:row-subtraction}
 \partial_{x_l}(Hf)(x)
 =\int\partial_{x_l}H(x,\widetilde x)
                         [f(\widetilde x)-f(x)]\dd\widetilde x.
\end{equation}
There is no derivative of $f(x)$, since this is an equality obtained
from the differentiated operator and the zero differentiated row.
One may first use smooth $f$ and then pass by continuity.

For $h_b\mathcal B_bg_b$ in \eqref{eq:response-kernels},
\begin{equation*}
\partial_{x_l}(h_b\mathcal B_bg_b)
 =h_b\left(D^2g_b[e_l,\delta u_b]+Dg_b\cdot\partial_lu_b(x)\right).
\end{equation*}
Since $h_b\norm{Du_b}_\infty\le C\kappa$ and
$h_b|\delta u_b|\le C\kappa\ell_b$, its kernel derivative is
bounded by
\begin{equation*}
C_\gamma\kappa\left(r^{-d+\beta-2\gamma-1}
           +\ell_br^{-d+\beta-2\gamma-2}\right)
                  \ind_{\{r\ge\ell_b/(2L)\}}.
\end{equation*}
Use \eqref{eq:row-subtraction}, Minkowski, and \eqref{eq:fractional-translation}. This gives
\begin{equation}\label{eq:H1-radial}
 \begin{aligned}
 &\norm{\nabla[(h_b\mathcal B_bg_b)f]}_2\\
 &\quad\le C_\gamma\kappa\norm f_{\dot H^\gamma}
 \left(\int_{\ell_b/(2L)<r<2r_0}r^{\beta-\gamma-2}\dd r
       +\ell_b\int_{\ell_b/(2L)<r<2r_0}r^{\beta-\gamma-3}\dd r\right)
 \le C_{\gamma,L}\kappa\ell_b^{-\alpha-\gamma}
                                             \norm f_{\dot H^\gamma}.
 \end{aligned}
\end{equation}
Here $\beta-\gamma-2=-1-\alpha-\gamma<-1$; integrating
from $\ell_b/(2L)$ to infinity bounds both expressions by
$C_{\gamma,L}\ell_b^{-\alpha-\gamma}$.
Therefore \eqref{eq:H1-radial} yields
\begin{equation}\label{eq:H1-summed-gradient}
 \norm{\nabla(H_1f)}_2 \le C_{\gamma,L}\kappa\sum_{b=0}^J \ell_b^{-\alpha-\gamma}\norm f_{\dot H^\gamma} \le C_{\gamma,L}\kappa \ell_J^{-\alpha-\gamma}\norm f_{\dot H^\gamma}.
\end{equation}

For $q_{ab}\mathcal B_a\mathcal B_bg_b$ in \eqref{eq:response-kernels}, differentiation gives
\begin{equation}\label{eq:H2-derivative}
 \begin{aligned}
 &\partial_{x_l}(q_{ab}\mathcal B_a\mathcal B_bg_b)
 =q_{ab}\bigl\{D^3g_b[e_l,\delta u_a,\delta u_b]
       +D^2g_b[\partial_lu_a(x),\delta u_b]\\
 &+D^2g_b[\delta u_a,\partial_lu_b(x)]
       +D^2g_b[e_l,\delta(\mathcal B_au_b)]
       +Dg_b\cdot\partial_l(\mathcal B_au_b)(x)\bigr\}.
 \end{aligned}
\end{equation}
If $a\ge b$, then $q_{ab}\le\kappa^2/(\lambda_a\lambda_b)$.
The bounds \eqref{eq:spatial}, \eqref{eq:velocity-increment},
\eqref{eq:material}, and \eqref{eq:material-gradient}, together
with $A_m/\lambda_m=\ell_m$, give the five estimates
\begin{equation*}
 q_{ab}|\delta u_a||\delta u_b|
 \le C\kappa^2\ell_a\ell_b,
\end{equation*}
\begin{equation*}
 q_{ab}\norm{Du_a}_\infty|\delta u_b|
 \le C\kappa^2\ell_b,
\end{equation*}
\begin{equation*}
 q_{ab}|\delta u_a|\norm{Du_b}_\infty
 \le C\kappa^2\ell_a,
\end{equation*}
\begin{equation*}
 q_{ab}|\delta(\mathcal B_au_b)|
 \le C\kappa^2\ell_a,
\end{equation*}
\begin{equation*}
 q_{ab}\norm{D_y(\mathcal B_au_b)}_\infty
 \le C\kappa^2.
\end{equation*}
Since $\ell_a\le\ell_b$, substitution of the five bounds above
in \eqref{eq:H2-derivative}, followed by
\eqref{eq:row-subtraction} and \eqref{eq:fractional-translation}, gives
\begin{equation}\label{eq:H2-radial}
 \begin{aligned}
 &\norm{\nabla[(q_{ab}\mathcal B_a\mathcal B_bg_b)f]}_2 \le C_\gamma\kappa^2\norm f_{\dot H^\gamma}
 \biggl(\int_{\ell_b/(2L)<r<2r_0}r^{\beta-\gamma-2}\dd r
 +\ell_b\int_{\ell_b/(2L)<r<2r_0}r^{\beta-\gamma-3}\dd r\\
 &+\ell_b^2\int_{\ell_b/(2L)<r<2r_0}r^{\beta-\gamma-4}\dd r\biggr) \le C_{\gamma,L}\kappa^2\ell_b^{-\alpha-\gamma}\norm f_{\dot H^\gamma}.
 \end{aligned}
\end{equation}
If $a<b$, then $q_{ab}\le\kappa^2/\lambda_b^2$.
Using $|\delta u_a|\le Cr\lambda_a$ and
$\norm{\mathcal B_au_b}_\infty\le C\lambda_aA_b$ gives
\begin{equation*}
 q_{ab}|\delta u_a||\delta u_b|
 \le C\kappa^2\frac{\lambda_a}{\lambda_b}r\ell_b,
\end{equation*}
\begin{equation*}
 q_{ab}\norm{Du_a}_\infty|\delta u_b|
 \le C\kappa^2\frac{\lambda_a}{\lambda_b}\ell_b,
\end{equation*}
\begin{equation*}
 q_{ab}|\delta u_a|\norm{Du_b}_\infty
 \le C\kappa^2\frac{\lambda_a}{\lambda_b}r,
\end{equation*}
\begin{equation*}
 q_{ab}|\delta(\mathcal B_au_b)|
 \le C\kappa^2\frac{\lambda_a}{\lambda_b}\ell_b,
\end{equation*}
\begin{equation*}
 q_{ab}\norm{D_y(\mathcal B_au_b)}_\infty
 \le C\kappa^2\frac{\lambda_a}{\lambda_b}.
\end{equation*}
Substituting the five bounds for $a<b$
in \eqref{eq:H2-derivative}, and then using
\eqref{eq:row-subtraction}, gives
\begin{equation}\label{eq:H2-coarser-pair-gradient}
 \begin{aligned}
 &\norm{\nabla[(q_{ab}\mathcal B_a\mathcal B_bg_b)f]}_2
 \le C_\gamma\kappa^2\frac{\lambda_a}{\lambda_b}
 \norm f_{\dot H^\gamma}
 \left(\int_{\ell_b/(2L)<r<2r_0}r^{\beta-\gamma-2}\dd r
  +\ell_b\int_{\ell_b/(2L)<r<2r_0}r^{\beta-\gamma-3}\dd r\right)\\
 &\le C_{\gamma,L}\kappa^2\frac{\lambda_a}{\lambda_b}
                 \ell_b^{-\alpha-\gamma}\norm f_{\dot H^\gamma}.
 \end{aligned}
\end{equation}
The index sums are uniformly bounded as follows:
\begin{equation*}
 \sum_{b=0}^J\sum_{a=b}^J\ell_b^{-\alpha-\gamma}\le\ell_J^{-\alpha-\gamma}\sum_{h\ge0}(h+1)2^{-h(\alpha+\gamma)},
\end{equation*}
\begin{equation*}
 \sum_{b=0}^J\ell_b^{-\alpha-\gamma}\sum_{a<b}\lambda_a/\lambda_b\le C\sum_{b=0}^J\ell_b^{-\alpha-\gamma}\le C_\gamma\ell_J^{-\alpha-\gamma}.
\end{equation*}
Equations~\eqref{eq:H2-radial} and
\eqref{eq:H2-coarser-pair-gradient} therefore give
\begin{equation}\label{eq:H2-summed-gradient}
 \norm{\nabla(H_2f)}_2
 \le C_{\gamma,L}\kappa^2\norm f_{\dot H^\gamma}
 \sum_{b=0}^J\ell_b^{-\alpha-\gamma}
       \left[J-b+1+\sum_{a<b}\frac{\lambda_a}{\lambda_b}\right] \le C_{\gamma,L}\kappa^2\ell_J^{-\alpha-\gamma}\norm f_{\dot H^\gamma}.
\end{equation}
All kernels used here
are regular at fixed $J,L$, so the differentiations and row
subtractions are justified as ordinary integrals.
Zero row and column integrals give, respectively,
\begin{equation*}
 H_i1=0,\qquad \int H_if=0.
\end{equation*}
Periodic Poincare, \eqref{eq:H1-summed-gradient} and
\eqref{eq:H2-summed-gradient} yield
\begin{equation*}
 \norm{H_if}_{H^1}
 \le C\norm{\nabla(H_if)}_2
 \le C_{\gamma,L}\kappa^i\ell_J^{-\alpha-\gamma}\norm f_{\dot H^\gamma}.
\end{equation*}
\end{proof}

\begin{lemma}\label{lem:stopped}
For every $\varepsilon>0$, $f\in H^1(\T^d)$ real, and $i=1,2$,
\begin{equation*}
 |\mu F_{\Delta_YH_i}(f)|+|F_{\mathcal B_vH_i}(f)|\le2\varepsilon\kappa\norm f_{\dot H^\gamma}^2+C_{\varepsilon,\gamma,L}\kappa\mu\norm{\nabla f}_2^2,
\end{equation*}
\begin{equation*}
 |F_{\mathcal B_vg}(f)|\le C_\gamma\mu\norm{\nabla f}_2^2.
\end{equation*}
The bilinear operators associated with $\Delta_YH_i$ and
$\mathcal B_vH_i$ have regular kernels; $F_{\mathcal B_vg}$ is
the Schur limit at fixed $N$ in Definition~\ref{def:responses}.
\end{lemma}
\begin{proof}
Symmetry and two spatial integrations by parts give
\begin{equation*}
\mu F_{\Delta_YH_i}(f)=-2\mu\langle\nabla f,\nabla(H_if)\rangle.
\end{equation*}
For a regular symmetric $H$, the operator $\mathcal B_v$
differentiates only $x,\widetilde x$. To estimate its quadratic form,
we use the skew potentials of \cite[Lemma 2.4, (2.5)]{FM}.
Subtracting the potentials for the scales above $J$ and above $N$
gives a matrix field $B$ such that
\begin{equation}\label{eq:potential}
 v=\diver B,\qquad B^T=-B,\qquad
 \norm B_\infty\le C\ell_J^{1+\alpha}\le C\mu,
\end{equation}
where $(\diver B)_i=\sum_m\partial_m B_{im}$.
Integration by parts gives
\begin{equation}\label{eq:skew}
 \begin{aligned}
 &F_{\mathcal B_vH}(f)
 =2\int f v_i\partial_i(Hf)\dd x =-2\int B_{im}\partial_m f\,\partial_i(Hf)\dd x
       -2\int f B_{im}\partial_{mi}(Hf)\dd x\\
 &=-2\int B_{im}\partial_m f\,\partial_i(Hf)\dd x.
 \end{aligned}
\end{equation}
In \eqref{eq:skew} the integral
$-2\int fB_{im}\partial_{mi}(Hf)\dd x$ vanishes by skew symmetry
of $B$ and symmetry of the Hessian. Thus each of
$|\mu F_{\Delta_YH_i}(f)|$ and $|F_{\mathcal B_vH_i}(f)|$
is at most
\begin{equation}\label{eq:H-pairing-bound}
C_{\gamma,L}\mu\kappa^i\ell_J^{-\alpha-\gamma}
                              \norm f_{\dot H^\gamma}\norm{\nabla f}_2.
\end{equation}
The choice of $J$ and $\beta>2\gamma$ imply
$\mu\ell_J^{-2(\alpha+\gamma)}\le C_\gamma$.
Young's inequality applied to \eqref{eq:H-pairing-bound} therefore gives
\begin{equation*}
 \begin{aligned}
 &|\mu F_{\Delta_YH_i}(f)|+|F_{\mathcal B_vH_i}(f)|
 \le2\varepsilon\kappa\norm f_{\dot H^\gamma}^2
 +C_{\varepsilon,\gamma,L}\kappa^{2i-1}\mu^2
 \ell_J^{-2(\alpha+\gamma)}\norm{\nabla f}_2^2\\
 &\le2\varepsilon\kappa\norm f_{\dot H^\gamma}^2
 +C_{\varepsilon,\gamma,L}\kappa\mu\norm{\nabla f}_2^2.
 \end{aligned}
\end{equation*}
For $H=g_\delta$, \eqref{eq:skew} holds and
$\norm{\nabla(g_\delta*f)}_2\le\norm{g_\delta}_1\norm{\nabla f}_2$.
Pass to $g$ by its $L^1$ convergence and the fixed-$N$ Schur convergence of $\mathcal B_vg_\delta$. Using \eqref{eq:potential}, this gives
\begin{equation*}
 |F_{\mathcal B_vg}(f)|
 \le2\norm B_\infty\norm g_1\norm{\nabla f}_2^2
 \le C_\gamma\mu\norm{\nabla f}_2^2.
\end{equation*}
\end{proof}

\section{Completion of the proof}\label{sec:completion}
We combine the preceding estimates in the following proposition,
then choose the parameters and pass to the infinite field and $L^2$ data.

\begin{proposition}[Finite-field estimate]\label{prop:finite-balance}
There are $\kappa_0>0$ and $c_\gamma,C_\gamma>0$ such that, for
$0<\kappa\le\kappa_0$, $L\ge4$, $0<\varepsilon<1$, $0<\mu\le1$,
$N\ge J(\mu)$, and a smooth deterministic datum, the scalar
$\theta^{\mu,N}$ driven by $u^{(N)}$ satisfies, with $T$ as in
\eqref{eq:model-horizon},
\begin{equation}\label{eq:finite-balance}
 \left\{c_\gamma\kappa
      -C_\gamma[L^{-\beta}+\kappa^2(1+\log L)^2]
      -4\varepsilon\kappa\right\}
    \E\int_0^T\norm{\theta^{\mu,N}_t}_{\dot H^\gamma}^2\dd t
       \le C_{\varepsilon,\gamma,L}\norm{\theta_0}_2^2.
\end{equation}
Here $c_\gamma,C_\gamma$ and $\kappa_0$ depend only on $\gamma$ and
the fixed model. The constant $C_{\varepsilon,\gamma,L}$ is independent
of $\kappa,\mu,N,T$ and the datum.
\end{proposition}

\begin{proof}[Proof of Proposition~\ref{prop:finite-balance}]
Choose $\kappa_0\le\min(\kappa_*,1)$ so that
Lemma~\ref{lem:coercivity} applies. Multiplication of
\eqref{eq:covariance-control} by $c_\gamma\kappa$ gives
\begin{equation}\label{eq:balance-multiplied-control}
 \begin{aligned}
 &c_\gamma\kappa\E\int_0^T\norm{\theta_t}_{\dot H^\gamma}^2\dd t
 \le-\mathrm{(C)}-\mu\E\int_0^TF_{\Delta_Yg}(\theta_t)\dd t\\
 &\quad+C_\gamma[\kappa L^{-\beta}+\kappa^2(1+\log L)]
                    \E\int_0^T\norm{\theta_t}_{\dot H^\gamma}^2\dd t\\
 &\quad+C_\gamma\E\int_0^T\norm{\theta_t}_2^2\dd t
       +C_\gamma\kappa^2(1+\log L)\mu
                              \E\int_0^T\norm{\nabla\theta_t}_2^2\dd t.
 \end{aligned}
\end{equation}
Solving \eqref{eq:renewal} for the two quantities on the right gives
\begin{equation*}
 \begin{aligned}
 &-\mathrm{(C)}-\mu\E\int_0^TF_{\Delta_Yg}(\theta_t)\dd t
 =F_g(\theta_0)-\E F_g(\theta_T)
       +\E\int_0^T[F_{E_0}(\theta_t)+F_{\mathcal B_vg}(\theta_t)]\dd t\\
 &+\E\int_0^T\sum_{a,b,c}q_{ab}(t)
             F_{\mathcal B_c\mathcal B_a\mathcal B_bg_b}(\theta_t)\dd t +\sum_{i=1}^2\left[\mu\E\int_0^TF_{\Delta_YH_i}(\theta_t)\dd t
                  +\E\int_0^TF_{\mathcal B_vH_i}(\theta_t)\dd t\right].
 \end{aligned}
\end{equation*}
Lemma~\ref{lem:transport} gives
\begin{equation*}
 \E\int_0^T\left[|F_{E_0}(\theta_t)|+
   \left|\sum_{a,b,c}q_{ab}(t)
         F_{\mathcal B_c\mathcal B_a\mathcal B_bg_b}(\theta_t)\right|\right]\dd t \le C_\gamma[L^{-\beta}+\kappa^2(1+\log L)^2]
                    \E\int_0^T\norm{\theta_t}_{\dot H^\gamma}^2\dd t.
\end{equation*}
Lemma~\ref{lem:stopped} gives the two bounds
\begin{equation*}
 \begin{aligned}
 &\sum_{i=1}^2\left[\mu\E\int_0^T|F_{\Delta_YH_i}(\theta_t)|\dd t
                  +\E\int_0^T|F_{\mathcal B_vH_i}(\theta_t)|\dd t\right]\\
 &\quad\le4\varepsilon\kappa
                    \E\int_0^T\norm{\theta_t}_{\dot H^\gamma}^2\dd t
       +C_{\varepsilon,\gamma,L}\kappa\mu
                    \E\int_0^T\norm{\nabla\theta_t}_2^2\dd t,
 \end{aligned}
\end{equation*}
\begin{equation*}
 \E\int_0^T|F_{\mathcal B_vg}(\theta_t)|\dd t
 \le C_\gamma\mu\E\int_0^T\norm{\nabla\theta_t}_2^2\dd t.
\end{equation*}
Combining these estimates with \eqref{eq:balance-multiplied-control},
and using $\kappa\le1$ and $L\ge4$, gives
\begin{equation}\label{eq:balance}
 \begin{aligned}
 &c_\gamma\kappa\E\int_0^T\norm{\theta_t}_{\dot H^\gamma}^2\dd t
 \le |F_g(\theta_0)|+\E|F_g(\theta_T)| +C_\gamma\E\int_0^T\norm{\theta_t}_2^2\dd t
      +C_{\varepsilon,\gamma,L}\mu
                            \E\int_0^T\norm{\nabla\theta_t}_2^2\dd t\\
 &+\left\{C_\gamma[L^{-\beta}+\kappa^2(1+\log L)^2]
                         +4\varepsilon\kappa\right\}
                        \E\int_0^T\norm{\theta_t}_{\dot H^\gamma}^2\dd t.
 \end{aligned}
\end{equation}
The endpoint terms are bounded by
$2\norm g_1\norm{\theta_0}_2^2$. The integrals of
$\norm{\theta_t}_2^2$ and $\mu\norm{\nabla\theta_t}_2^2$ are
controlled by \eqref{eq:energy}, since $T\le2$.
Moving the Sobolev integral on the right of \eqref{eq:balance}
to the left gives
\eqref{eq:finite-balance}.
\end{proof}

\begin{proof}[Proof of Theorem~\ref{thm:main}]
Fix $\varepsilon>0$ sufficiently small depending on $\gamma$, and set
\begin{equation*}
 L=\kappa^{-2/\beta}.
\end{equation*}
Then $L^{-\beta}+\kappa^2(1+\log L)^2=o(\kappa)$ as
$\kappa\downarrow0$. We can therefore choose
$\kappa_\gamma\le\kappa_0$ so that $L\ge4$, all preceding
smallness conditions hold, and
\begin{equation}\label{eq:absorption-smallness}
 C_\gamma[L^{-\beta}+\kappa^2(1+\log L)^2]+4\varepsilon\kappa
 \le\tfrac12c_\gamma\kappa.
\end{equation}
Proposition~\ref{prop:finite-balance} and \eqref{eq:absorption-smallness} give, for all
$N\ge J(\mu)$,
\begin{equation}\label{eq:finite}
 \E\int_0^1\norm{\theta^{\mu,N}_t}_{\dot H^\gamma}^2\dd t
 \le\E\int_0^T\norm{\theta^{\mu,N}_t}_{\dot H^\gamma}^2\dd t
       \le C_{\gamma,\kappa}\norm{\theta_0}_2^2.
\end{equation}
The comparison in \eqref{eq:finite} uses nonnegativity and $T\ge1$.
For $0<\gamma<1$, the Fourier weights give
\begin{equation*}
 \norm f_{H^\gamma}^2\le\norm f_2^2+\norm f_{\dot H^\gamma}^2.
\end{equation*}
Restricting to $[0,1]$ and using \eqref{eq:energy}, we obtain
\begin{equation}\label{eq:finite-inhomogeneous}
 \E\int_0^1\norm{\theta_t^{\mu,N}}_{H^\gamma}^2\dd t
 \le(1+C_{\gamma,\kappa})\norm{\theta_0}_2^2.
\end{equation}
The restriction solves \eqref{eq:PDE-2} on $[0,1]$. Constants may
depend on the fixed $L$ and hence on $\kappa$, but not on $N,\mu$.

We complete the proof of the theorem by passing to the limit
as $N\to\infty$ in \eqref{eq:finite}. If $w,\widetilde w$
are bounded divergence-free drifts with the same datum and solutions
$\theta,\widetilde\theta$, their difference $h=\theta-\widetilde\theta$
satisfies
\begin{equation*}
\partial_t h+w\cdot\nabla h=\mu\Delta h
                  -\diver((w-\widetilde w)\widetilde\theta),\qquad h(0)=0.
\end{equation*}
Testing with $h$ and applying Young's inequality gives
\begin{equation*}
 \tfrac12\frac\dd{\dd t}\norm h_2^2+\mu\norm{\nabla h}_2^2 =\langle(w-\widetilde w)\widetilde\theta,\nabla h\rangle \le\tfrac\mu2\norm{\nabla h}_2^2 +\frac{\norm{w-\widetilde w}_\infty^2}{2\mu} \norm{\widetilde\theta}_2^2.
\end{equation*}
Integrating over $[0,t]\subset[0,1]$ and using
\eqref{eq:energy} for $\widetilde\theta$ gives
\begin{equation}\label{eq:stability}
 \frac12\norm{h(t)}_2^2+\frac\mu2\int_0^t\norm{\nabla h}_2^2\dd s \le\frac{\norm{w-\widetilde w}_\infty^2}{2\mu} \int_0^t\norm{\widetilde\theta}_2^2\dd s \le\frac{\norm{w-\widetilde w}_\infty^2}{2\mu}\norm{\theta_0}_2^2.
\end{equation}
The bound \eqref{eq:spatial} gives
\begin{equation*}
 \norm{u^{(N)}-u^{(\infty)}}_\infty
 \le C\sum_{b>N}2^{-\alpha b}
 \le C2^{-\alpha N}.
\end{equation*}
Thus \eqref{eq:stability}, with $w=u^{(N)}$ and
$\widetilde w=u^{(\infty)}$, yields, at fixed $\mu>0$,
\begin{equation*}
 \E\sup_{0\le t\le1}\norm{\theta_t^{\mu,N}-\theta_t^\mu}_2^2
 \le\frac1\mu\E\norm{u^{(N)}-u^{(\infty)}}_\infty^2\norm{\theta_0}_2^2
 \le\frac{C2^{-2\alpha N}}\mu\norm{\theta_0}_2^2\longrightarrow0.
\end{equation*}
The same deterministic bound in \eqref{eq:stability} also gives
pathwise convergence in $C([0,1];L^2)$. Each finite-field solution
is measurable by \cite[Lemma 9.1]{FM}; its limit $\theta^\mu$ is
therefore measurable as well.
For a fixed Fourier radius $Q$,
\begin{equation*}
 \E\int_0^1\sum_{|n|\le Q}(1+4\pi^2|n|^2)^\gamma
     |\widehat{\theta^{\mu,N}}(n,t)-\widehat{\theta^\mu}(n,t)|^2\dd t \le(1+4\pi^2Q^2)^\gamma
     \E\sup_{0\le t\le1}\norm{\theta_t^{\mu,N}-\theta_t^\mu}_2^2
 \longrightarrow0.
\end{equation*}
By Cauchy--Schwarz and \eqref{eq:energy},
\begin{equation}\label{eq:Fourier-energy-convergence}
 \begin{aligned}
 &\left|\E\int_0^1\sum_{|n|\le Q}(1+4\pi^2|n|^2)^\gamma
     (|\widehat{\theta^{\mu,N}}(n,t)|^2-|\widehat{\theta^\mu}(n,t)|^2)\dd t\right|\\
 &\quad\le2(1+4\pi^2Q^2)^\gamma\norm{\theta_0}_2
     \left(\E\sup_{0\le t\le1}\norm{\theta_t^{\mu,N}-\theta_t^\mu}_2^2\right)^{1/2}
 \longrightarrow0.
 \end{aligned}
\end{equation}
Consequently, \eqref{eq:finite-inhomogeneous} gives

\begin{equation*}
 \begin{aligned}
 &\E\int_0^1\sum_{|n|\le Q}(1+4\pi^2|n|^2)^\gamma
                         |\widehat{\theta^\mu}(n,t)|^2\dd t\\
 &\qquad=\lim_{N\to\infty}\E\int_0^1\sum_{|n|\le Q}
 (1+4\pi^2|n|^2)^\gamma|\widehat{\theta^{\mu,N}}(n,t)|^2\dd t
 \le C_{\gamma,\kappa}\norm{\theta_0}_2^2.
 \end{aligned}
\end{equation*}
Monotone convergence as $Q\to\infty$ transfers the full estimate.
Finally choose smooth data $\theta_{0,h}\to\theta_0$ in $L^2$.
Contraction gives
\begin{equation*}
 \sup_{0\le t\le1}\norm{\theta_h^\mu(t)-\theta^\mu(t)}_2
 \le\norm{\theta_{0,h}-\theta_0}_2.
\end{equation*}
The same computation as in \eqref{eq:Fourier-energy-convergence} yields
\begin{equation*}
 \begin{aligned}
 &\left|\E\int_0^1\sum_{|n|\le Q}(1+4\pi^2|n|^2)^\gamma
       (|\widehat{\theta_h^\mu}(n,t)|^2-|\widehat{\theta^\mu}(n,t)|^2)\dd t\right|\\
 &\quad\le(1+4\pi^2Q^2)^\gamma
      (\norm{\theta_{0,h}}_2+\norm{\theta_0}_2)
                           \norm{\theta_{0,h}-\theta_0}_2\longrightarrow0.
 \end{aligned}
\end{equation*}
Passing to the limit in the finite Fourier sum therefore gives
\begin{equation*}
 \E\int_0^1\sum_{|n|\le Q}(1+4\pi^2|n|^2)^\gamma
                         |\widehat{\theta^\mu}(n,t)|^2\dd t
 \le C_{\gamma,\kappa}\norm{\theta_0}_2^2.
\end{equation*}
Monotone convergence as $Q\to\infty$ gives
\begin{equation*}
 \E\int_0^1\norm{\theta_t^\mu}_{H^\gamma}^2\dd t
 \le C_{\gamma,\kappa}\norm{\theta_0}_2^2,\qquad0<\mu\le1.
\end{equation*}
\end{proof}

\bibliographystyle{amsplain}

\begin{thebibliography}{10}

\bibitem{ACMMixing}
Giovanni Alberti, Gianluca Crippa, and Anna~L. Mazzucato, \emph{Exponential
  self-similar mixing by incompressible flows}, Journal of the American
  Mathematical Society \textbf{32} (2019), no.~2, 445--490.

\bibitem{ACMLoss}
\bysame, \emph{Loss of regularity for the continuity equation with
  non-{Lipschitz} velocity field}, Annals of PDE \textbf{5} (2019), no.~1, 9.

\bibitem{AV}
Scott Armstrong and Vlad Vicol, \emph{Anomalous diffusion by fractal
  homogenization}, Annals of PDE \textbf{11} (2025), 2.

\bibitem{AVSurvey}
\bysame, \emph{Anomalous diffusion via iterative quantitative homogenization:
  an overview of the main ideas}, 2025,
  \href{https://arxiv.org/abs/2503.11744}{arXiv:2503.11744}.

\bibitem{BBPSBatchelor}
Jacob Bedrossian, Alex Blumenthal, and Samuel Punshon-Smith, \emph{The
  {Batchelor} spectrum of passive scalar turbulence in stochastic fluid
  mechanics at fixed {Reynolds} number}, Communications on Pure and Applied
  Mathematics \textbf{75} (2022), no.~6, 1237--1291.

\bibitem{BBPSLagrangian}
\bysame, \emph{{Lagrangian} chaos and scalar advection in stochastic fluid
  mechanics}, Journal of the European Mathematical Society \textbf{24} (2022),
  no.~6, 1893--1990.

\bibitem{BGK}
Denis Bernard, Krzysztof Gaw{\c e}dzki, and Antti Kupiainen, \emph{Slow modes
  in passive advection}, Journal of Statistical Physics \textbf{90} (1998),
  no.~3--4, 519--569.

\bibitem{BCCDL}
Elia Bru{\`e}, Maria Colombo, Gianluca Crippa, Camillo {De Lellis}, and Massimo
  Sorella, \emph{{Onsager} critical solutions of the forced {Navier--Stokes}
  equations}, Communications on Pure and Applied Analysis \textbf{23} (2024),
  no.~10, 1350--1366.

\bibitem{BDL}
Elia Bru{\`e} and Camillo {De Lellis}, \emph{Anomalous dissipation for the
  forced {3D Navier--Stokes} equations}, Communications in Mathematical Physics
  \textbf{400} (2023), no.~3, 1507--1533.

\bibitem{BSWRegularity}
Jan Burczak, L{\'a}szl{\'o} Sz{\'e}kelyhidi, and Bian Wu, \emph{Scalar
  anomalous dissipation and optimal regularity via iterated homogenization},
  2026, \href{https://arxiv.org/abs/2604.13912}{arXiv:2604.13912}.

\bibitem{BSWEuler}
Jan Burczak, L{\'a}szl{\'o} Sz{\'e}kelyhidi, Jr., and Bian Wu, \emph{Anomalous
  dissipation and {Euler} flows}, 2023,
  \href{https://arxiv.org/abs/2310.02934}{arXiv:2310.02934}.

\bibitem{CCS}
Maria Colombo, Gianluca Crippa, and Massimo Sorella, \emph{Anomalous
  dissipation and lack of selection in the {Obukhov--Corrsin} theory of scalar
  turbulence}, Annals of PDE \textbf{9} (2023), no.~2, 21.

\bibitem{CHR}
Maria Colombo, Elias Hess-Childs, and Keefer Rowan, \emph{Effective
  {Lagrangian} regularity and the uniqueness threshold for random {H{\"o}lder}
  velocity fields}, 2026,
  \href{https://arxiv.org/abs/2608.05931}{arXiv:2608.05931}.

\bibitem{CKRZ}
Peter Constantin, Alexander Kiselev, Lenya Ryzhik, and Andrej Zlato{\v s},
  \emph{Diffusion and mixing in fluid flow}, Annals of Mathematics \textbf{168}
  (2008), no.~2, 643--674.

\bibitem{Corrsin}
Stanley Corrsin, \emph{On the spectrum of isotropic temperature fluctuations in
  an isotropic turbulence}, Journal of Applied Physics \textbf{22} (1951),
  no.~4, 469--473.

\bibitem{DEIJ}
Theodore~D. Drivas, Tarek~M. Elgindi, Gautam Iyer, and In-Jee Jeong,
  \emph{Anomalous dissipation in passive scalar transport}, Archive for
  Rational Mechanics and Analysis \textbf{243} (2022), no.~3, 1151--1180.

\bibitem{DE}
Theodore~D. Drivas and Gregory~L. Eyink, \emph{A {Lagrangian}
  fluctuation--dissipation relation for scalar turbulence. {Part I}. {Flows}
  with no bounding walls}, Journal of Fluid Mechanics \textbf{829} (2017),
  153--189.

\bibitem{DGP}
Theodore~D. Drivas, Lucio Galeati, and Umberto Pappalettera, \emph{Anomalous
  dissipation and regularization in isotropic {Gaussian} turbulence}, 2025,
  \href{https://arxiv.org/abs/2509.10211v1}{arXiv:2509.10211v1}.

\bibitem{EL}
Tarek~M. Elgindi and Kyle Liss, \emph{Norm growth, non-uniqueness, and
  anomalous dissipation in passive scalars}, Archive for Rational Mechanics and
  Analysis \textbf{248} (2024), no.~6, 120.

\bibitem{FGV}
Gregory Falkovich, Krzysztof Gaw{\c e}dzki, and Massimo Vergassola,
  \emph{Particles and fields in fluid turbulence}, Reviews of Modern Physics
  \textbf{73} (2001), no.~4, 913--975.

\bibitem{FM}
Charles Fefferman and Svitlana Mayboroda, \emph{Richardson dispersion and
  anomalous dissipation for the random velocity model with controlled
  sweeping}, Forthcoming preprint.

\bibitem{GGM}
Lucio Galeati, Francesco Grotto, and Mario Maurelli, \emph{Anomalous
  regularization in {Kraichnan}'s passive scalar model}, Probability Theory and
  Related Fields \textbf{196} (2026), 255--291.

\bibitem{GK}
Krzysztof Gaw{\c e}dzki and Antti Kupiainen, \emph{Anomalous scaling of the
  passive scalar}, Physical Review Letters \textbf{75} (1995), no.~21,
  3834--3837.

\bibitem{Heil}
Christopher Heil, \emph{Metrics, norms, inner products and operator theory:
  {Chapter} 8, integral operators}, 2020, Online extra chapter, updated October
  3, 2020, Section 8.3, Theorems 8.3.1 and 8.3.4(b), pp. 363--365.
  \url{https://heil.math.gatech.edu/metricnote/chap8.pdf}.

\bibitem{HCRIntermittent}
Elias Hess-Childs and Keefer Rowan, \emph{Turbulent and intermittent phenomena
  in a universal total anomalous dissipator}, 2025,
  \href{https://arxiv.org/abs/2508.00115}{arXiv:2508.00115}.

\bibitem{HCRUniversal}
\bysame, \emph{A universal total anomalous dissipator}, Communications in
  Mathematical Physics \textbf{407} (2026), 257.

\bibitem{HPZZ}
Martina Hofmanov{\'a}, Umberto Pappalettera, Rongchan Zhu, and Xiangchan Zhu,
  \emph{Anomalous and total dissipation due to advection by solutions of
  randomly forced {Navier--Stokes} equations}, The Annals of Applied
  Probability \textbf{35} (2025), no.~5, 3119--3149.

\bibitem{JDEMS}
R.~Jimenez, S.~Mayboroda, C.~{De Lellis}, D.~Spergel, and M.~Engelstein,
  \emph{Eddies that ride the flow: a model of turbulence with {Richardson}
  dispersion and anomalous dissipation}, Forthcoming preprint.

\bibitem{JSAutonomous}
Carl Johan~Peter Johansson and Massimo Sorella, \emph{Anomalous dissipation via
  spontaneous stochasticity with a two-dimensional autonomous velocity field},
  Duke Mathematical Journal \textbf{175} (2026), no.~11,
  \href{https://doi.org/10.1215/00127094-2025-0067}{doi:10.1215/00127094-2025-0067}.

\bibitem{JSMeasure}
\bysame, \emph{Nontrivial absolutely continuous part of anomalous dissipation
  measures in time}, Journal of Differential Equations \textbf{453} (2026),
  113912.

\bibitem{Kraichnan}
Robert~H. Kraichnan, \emph{Small-scale structure of a scalar field convected by
  turbulence}, The Physics of Fluids \textbf{11} (1968), no.~5, 945--953.

\bibitem{LJR}
Yves {Le Jan} and Olivier Raimond, \emph{Integration of {Brownian} vector
  fields}, The Annals of Probability \textbf{30} (2002), no.~2, 826--873.

\bibitem{Obukhov}
A.~M. Obukhov, \emph{The structure of the temperature field in a turbulent
  flow}, Izvestiya Akademii Nauk SSSR, Seriya Geograficheskaya i
  Geofizicheskaya \textbf{13} (1949), no.~1, 58--69, In Russian.

\bibitem{RS}
Luz Roncal and Pablo~Ra{\'u}l Stinga, \emph{Fractional {Laplacian} on the
  torus}, Communications in Contemporary Mathematics \textbf{18} (2016),
  1550033, Theorem numbering refers to
  \href{https://arxiv.org/abs/1209.6104v3}{arXiv:1209.6104v3}.

\bibitem{RowanDiffusion}
Keefer Rowan, \emph{On anomalous diffusion in the {Kraichnan} model and
  correlated-in-time variants}, Archive for Rational Mechanics and Analysis
  \textbf{248} (2024), no.~5, 93.

\bibitem{RowanSpectrum}
\bysame, \emph{The {Obukhov--Corrsin} spectrum of passive scalar turbulence
  through anomalous regularization}, 2025,
  \href{https://arxiv.org/abs/2512.02853}{arXiv:2512.02853}.

\bibitem{ShraimanSiggia}
Boris~I. Shraiman and Eric~D. Siggia, \emph{Scalar turbulence}, Nature
  \textbf{405} (2000), no.~6787, 639--646.

\bibitem{Sreenivasan}
Katepalli~R. Sreenivasan, \emph{Turbulent mixing: a perspective}, Proceedings
  of the National Academy of Sciences of the United States of America
  \textbf{116} (2019), no.~37, 18175--18183.

\bibitem{Warhaft}
Zellman Warhaft, \emph{Passive scalars in turbulent flows}, Annual Review of
  Fluid Mechanics \textbf{32} (2000), 203--240.

\bibitem{Yaglom}
A.~M. Yaglom, \emph{On the local structure of a temperature field in a
  turbulent flow}, Doklady Akademii Nauk SSSR \textbf{69} (1949), no.~6,
  743--746, In Russian.

\end{thebibliography}
\providecommand{\bysame}{\leavevmode\hbox to3em{\hrulefill}\thinspace}
\providecommand{\MR}{\relax\ifhmode\unskip\space\fi MR }
\providecommand{\MRhref}[2]{%
  \href{http://www.ams.org/mathscinet-getitem?mr=#1}{#2}
}
\providecommand{\href}[2]{#2}

\end{document}